\magnification=1200

\def\A{{\cal A}}

\def\D{{\cal D}}
\def\E{{\cal E}}

\def\N{{\bf N}}
\def\Q{{\bf Q}}
\def\R{{\bf R}}
\def\X{{\cal X}}

\def\hal{{\vrule height 10pt width 4pt depth 0pt}}
\def\la{{\langle}}
\def\ra{{\rangle}}
\def\Li{{\rm Lip}}

\def\O1{{\Omega^1}}

\centerline{\bf Lipschitz algebras and derivations II:
exterior differentiation}
\medskip

\centerline{Nik Weaver\footnote{}{
{\it Key words and phrases:}
Banach manifold, Banach module, derivation, Dirichlet space,
exterior derivative, Lipschitz algebra, metric space, rectifiable set,
sub-Riemannian manifold, vector field.\hfill\break
}}
\bigskip
\bigskip

{\narrower{
\noindent Basic aspects of differential geometry can be extended to various
non-classical settings: Lipschitz manifolds, rectifiable sets, sub-Riemannian
manifolds, Banach manifolds, Wiener space, etc. Although the constructions
differ, in each of these cases
one can define a module of measurable 1-forms and a
first-order exterior derivative. We give a general construction which
applies to any metric space equipped with a $\sigma$-finite measure and
produces the desired result in all of the above cases. It also applies to an
important class of Dirichlet spaces, where, however,
the known first-order differential calculus in
general differs from ours (although the two are related).
\bigskip}}
\bigskip

\noindent {\underbar{\bf 1. Introduction.}}
\bigskip

This paper is a continuation of [69]. There we considered derivations of
$L^\infty(X,\mu)$ into a certain kind of bimodule and constructed an
associated metric on $X$. Here we take a metric space $M$ as given and
consider only
derivations into modules whose left and right actions coincide (monomodules).
This allows the extraction of a kind of differentiable structure on $M$.
The exterior derivative is the universal derivation, into a certain type
of bimodule, which is compatible with the metric in a special sense.
\medskip

Surprisingly, this construction requires no serious conditions
on $M$; in particular, it need not be a manifold. But actually several
lines of research point in this direction. First, in Connes' noncommutative
geometry ([12], [13]) it appears that however one makes sense of the
noncommutative version of the property of
being a manifold, it is basically unrelated
to the noncommutative version of differentiable structure (and the
latter is, from a functional-analytic point of view, simpler). This algebraic
approach to differentiable structure is already explicit in
[57] and [58]. Second, Gromov ([29], [30]) has shown that
relatively sophisticated geometric notions, such as sectional curvature,
can be treated in a purely metric fashion (see also [9] and [56]).
Third, Sauvageot was able to construct an exterior derivative from initial
data consisting only of a Dirichlet
form with certain properties ([54], [55]), and others have also
treated Dirichlet spaces geometrically ([5], [60-63]). Finally,
recent work on fractals ([36], [37], [41], [53]) has a strong geometric flavor,
again indicating that manifold structure is not necessary for some kind of
differential analysis.
\medskip

Our main result is the identification of a module of measurable 1-forms
and a first-order exterior derivative
on any metric space equipped with a $\sigma$-finite
Borel measure. The consequences are availability of basic
differential geometric tools in non-classical settings, and
unification of previous work along these lines.
\medskip

Actually, the measure itself is not relevant to our construction, only its
measure class. In most finite-dimensional examples
Hausdorff measure (see e.g.\ [25]) provides the appropriate measure
class. To be precise, there is at most one value of $p \in [0,\infty)$
such that $p$-dimensional Hausdorff measure is $\sigma$-finite but nonzero,
and this is the desired measure. So in a restricted setting
one can regard our construction
as proceeding from the metric alone. However, we will want to consider some
infinite-dimensional examples (meaning that $p$-dimensional Hausdorff measure
is infinite for all $p$), and here it appears that the measure must be
added as an independent ingredient. There may also be occasional
finite-dimensional examples where one wants to use something other than
Hausdorff measure class (e.g.\ see the comment following Theorem 40 and the
second example considered in Theorem 55).
\medskip

It is interesting to note that for much of what we do the infinite-dimensional
Riemannian/Hilbert module case is just as tractable as the general
finite-dimensional case. (See Corollary 24 and the comment following
Theorem 10.)
\bigskip

\noindent {\bf A. Non-classical spaces}
\bigskip

By the ``classical'' setting we mean the case of a smooth (or at worst
$C^1$) Riemannian manifold. Applying our construction in this case,
using no structure besides the distance function and the Lebesgue measure
class, produces the usual (measurable) differential geometric data
associated with the manifold.
\medskip

More generally, let $M$ be a Lipschitz manifold. By this we mean
that $M$ is a metric space which is locally bi-Lipschitz equivalent to the
unit ball in $\R^n$; see ([17], [43], [52]) for background. These spaces still
possess measurable 1-forms and a first-order exterior derivative, and our
construction still recovers them.
\medskip

In a different direction, one can generalize Riemannian manifolds by
sub-Riemannian (or ``Carnot-Carath\'eodory'') spaces [4]. Here the metric
is defined in terms of paths which are tangent to a fixed subbundle of the
tangent bundle, and our construction then recovers this subbundle. (Cf.\
the more sophisticated view of tangent spaces for left-invariant sub-Riemannian
metrics on Lie groups described in [3] and [48].)
\medskip

Any sub-Riemannian manifold may be viewed as a Gromov-Hausdorff limit of
Riemannian manifolds [30]; but in general these limits need not
themselves be Riemannian, sub-Riemannian, or even manifolds. Nonetheless,
Peter Petersen has pointed out to me that one can define a measure on any such
limit by treating Gromov-Hausdorff convergence as Hausdorff convergence
within a larger space (see [50]) and taking a weak* limit of measures on the
approximating spaces. So our construction will also apply in this
situation, but I have not pursued this.
\medskip

Rectifiable sets in the sense of geometric measure theory ([22], [46])
constitute another broad class, similar to the Lipschitz manifolds,
in which a geometric approach has proved
valuable. Here too there is a natural tangent space defined almost everywhere
which can be recovered from the metric alone.
\medskip

Generalizing to infinite dimensions, one can treat metric spaces that
are locally bi-Lipschitz equivalent to the unit ball of a Banach space
and are equipped with a measure that satisfies a mild condition. Examples
of such spaces are path spaces of the type considered in [20] and
[47], for example. (But not spaces of continuous paths; see the discussion
of Weiner space below.)
\medskip

Non-rectifiable fractals --- specifically, fractals with non-integral
Hausdorff dimension --- have also been investigated in various ways that
have some geometric flavor ([34],
[36], [37], [38], [41], [64]). This is seen most
explicitly in the renormalized Laplace operator and Gauss-Green's formula
of [36]. There has also been interesting work on diffusion processes in
this setting ([1], [2], [39], [53]). Unfortunately, our construction is
generally vacuous for sets of this type, so it seems that we have nothing to
contribute here. I believe there is no meaningful ``renormalized''
tangent bundle for these sets.
\medskip

But the diffusion processes just mentioned lack an important property: the
associated Dirichlet forms do not admit a {\it carr\'e du champ} (see [8]).
In the presence of this extra hypothesis there is a natural underlying
metric and, independently, an elegant construction of an exterior derivative
[54]. But the latter in general is different from our exterior derivative.
Their relationship is addressed in Proposition 54.
\medskip

As a special case, one may consider the standard diffusion process and
associated Dirichlet form on Wiener space. This does admit a {\it carr\'e du
champ}, and the two exterior derivatives agree. This example has been
considered in detail and can be described in an elementary fashion, without
reference to Dirichlet forms (see [8], [66]). It is important to realize
that although Weiner space is a Banach manifold, this is not really relevant
to its differentiable structure: the Gross-Sobolev derivative is not directly
derived from the metric which provides
local Banach space structure. In particular,
the tangent spaces carry inner products despite the fact that Weiner space,
in its usual formulation, is not locally isomorphic to a Hilbert space.
\medskip

We must also mention
recent progress which has been made in several non-classical directions in
geometry which do not clearly fit with the point of view taken here. One of
these is Harrison's work on nonsmooth chains [32], which deals with
$p$-forms on nonsmooth sets. This is really a generalization in a different
direction and possibly could be combined with our approach.
Davies' analysis on graphs
([15], [16]), based on Connes' two-point space ([13], Example VI.3.1),
involves a notion of 1-form, but this is a discrete, not a
differential, object. (The definition of $p$-forms for the product of a
two-point space by an ordinary manifold is treated nicely in [40].)
Finally, noncommutative geometry in the sense of Connes ([12], [13])
has become a major industry, with applications in many directions.
It is unclear to what extent the idea of vector fields as derivations,
which is central to this paper, is helpful in the noncommutative setting. The
main reason for thinking not is that if $A$ is a noncommutative algebra, then
the set of derivations from $A$ into itself generally is not a module over
$A$. Still, some progress in this direction has been made (see [18], [45]).
\bigskip

\noindent {\bf B. Plan of the paper and acknowledgements}
\bigskip

We begin with some background material on the type of
$L^\infty(X)$-modules of interest to us. These ``abelian W*-modules'' are
introduced and their structure analyzed in section 2. Then in section 3 we
describe the notion of a ``measurable metric'' and its relation to derivations
into bimodules. These two sections set up the fundamental notions that
are needed in the sequel.
\medskip

The two succeeding sections are the core of the paper. Our construction
and its general properties are given in section 4. Section 5 is
devoted to examples; here we show that our construction produces the
standard result in most if not all of the cases where there is one, and we
also consider a few new cases.
Finally, in section 6 we carefully consider the case of Dirichlet spaces.
\medskip

The original source of motivation for this work was Sauvageot's pair of
papers [54] and [55], which show that one can have a meaningful exterior
derivative in the absence of anything resembling manifold structure.
It has also been my good fortune to have been given in person a great
deal of information and advice on the topics considered here. In this
regard Daniel Allcock, Renato Feres, Ken Goodearl, Gary Jensen,
Mohan Kumar, Michel Lapidus, Paul Muhly, Peter Petersen,
J\"urgen Schweizer, Mitch Taibleson,
and Ed Wilson all contributed to the mathematical content of this paper, and
it is a pleasure to acknowledge their help.
\medskip

My primary debt, however, is owed to Marty Silverstein, in whose seminar
I learned most of what I know about Dirichlet forms, and who encouraged
this project when it was in an early stage. My sincere thanks also go to him.
\bigskip
\bigskip

\noindent {\underbar{\bf 2. Abelian W*-modules.}}
\bigskip

\noindent {\bf A. Definitions}
\bigskip

The scalar field will be real throughout. We will invoke
several facts from the literature which involve complex scalars, but this
raises no serious issues. In every case a trivial complexification argument
justifies the application.
\medskip

Let $K$ be a compact Hausdorff space and let $E$ be a module over $C(K)$.
Recall that $E$ is a {\it Banach module} if it is also a Banach space and its
norm satisfies $\|f\phi\| \leq \|f\|\|\phi\|$ for all $f \in C(K)$
and $\phi \in E$.
\medskip

Following [19], we say that $E$ is a $C(K)$-normed module
if there exists a map $|\cdot|: E \to C(K)$ such that
$|\phi| \geq 0$, $\|\phi\| = \|\,|\phi|\,\|_\infty$, and $|f\phi|
= |f| |\phi|$ for all $f \in C(K)$ and $\phi \in E$. This map is to be
thought of as a fiberwise norm; according to ([51],
Corollary 6), it is unique if it exists. By ([19], p.\ 48) and
([51], Proposition 2) we have the following equivalent characterization of
$C(K)$-normed modules:
\bigskip

\noindent {\bf Theorem 1.} {\it Let $E$ be a Banach module over $C(K)$.
Then $E$ is a $C(K)$-normed module if and only if it is isometrically
isomorphic to the set of continuous sections of some (F) Banach bundle
over $K$.}\hfill\hal
\bigskip

Here, as in [19], (F) Banach bundles are Fell bundles $B$ of Banach spaces
$B_x$, i.e.\ $B = \bigcup_{x \in K}B_x$. Their defining property is that the
Banach space norm is continuous as a map from $B$ to $\R$.
\medskip

Now let $X = (X,\mu)$ be a measure space. We will assume throughout that
$X$ is $\sigma$-finite, although everything we do should work in
the more general case that $X$ is finitely decomposable. By this we mean
that $X$ can be expressed as a (possibly uncountable)
disjoint union of finite measure subsets, $X = \bigcup X_i$, such that
$A \subset X$ is measurable if and only if $A \cap X_i$ is measurable for
all $i$, in which case $\mu(A) = \sum \mu(A \cap X_i)$. (The spaces
$L^\infty(X)$ for $X$ finitely decomposable are precisely the real parts
of abelian von Neumann algebras.)
\medskip

Now $L^\infty(X)$ is isomorphic to $C(K)$ for some extremely
disconnected compact Hausdorff space $K$, so the notion of a $C(K)$-normed
module specializes to this case. In this situation we have a further
equivalence ([51], Theorem 9):
\bigskip

\noindent {\bf Theorem 2.} {\it Let $E$ be a Banach module over
$L^\infty(X)$. The following are equivalent:
\medskip

{\narrower{
\noindent (a). $E$ is an $L^\infty(X)$-normed module;
\medskip

\noindent (b). there is a compact Hausdorff space $K$, an isometric algebra
homomorphism $\imath$ from $L^\infty(X)$ into $C(K)$, and an isometric linear
map $\pi$ from $E$ into $C(K)$ such that $\pi(f\phi) = \imath(f)\pi(\phi)$
for all $f \in L^\infty(X)$ and $\phi \in E$;
\medskip}}

(c). $\|\phi\| = {\rm max}(\|p\phi\|, \|(1 - p)\phi\|)$ for all
$\phi \in E$ and any projection $p \in L^\infty(X)$.}\hfill\hal
\bigskip

Theorem 2 (b) is an abelian version of the notion of an operator module
(see [21]). Modules with this property were used heavily in [69] and are
central to this paper as well. The precise class of modules which we need
is specified in the next definition.
\bigskip

\noindent {\bf Definition 3.} A Banach module $E$ over $L^\infty(X)$ is an
{\it abelian W*-module} if it satisfies the equivalent conditions of
Theorem 2, and in addition is a dual Banach space such that the product map
$L^\infty(X)\times E \to E$ is separately weak*-continuous in each
variable.\hfill\hal
\bigskip

The notion of duality is central to our analysis of abelian W*-modules. This
concept is given in the next definition; following it, we give a module
version of the Hahn-Banach theorem.
\bigskip

\noindent {\bf Definition 4.} Let $E$ be a Banach module over $C(K)$.
Then its {\it dual module} $E'$ is the set of norm-bounded $C(K)$-module
homomorphisms from $E$ into $C(K)$. Give each such homomorphism its norm
as an operator between Banach spaces, and define the module operation by
$$(f\Phi)(\phi) = f\cdot \Phi(\phi)$$
for $f \in C(K)$, $\phi \in E$, and $\Phi \in E'$. It is easy to check that
$E'$ is a Banach module over $C(K)$.\hfill\hal
\bigskip

\noindent {\bf Theorem 5.} {\it Let $E$ be an $L^\infty(X)$-normed module,
let $E_0 \subset E$ be a submodule, and let $\Phi_0: E_0 \to L^\infty(X)$
be a norm-one module homomorphism. Then there is a norm-one module homomorphism
$\Phi: E \to L^\infty(X)$ such that $\Phi|_{E_0} = \Phi_0$.}
\medskip

\noindent {\it Proof.} By Zorn's lemma, it will suffice to consider the case
that there exists $\phi_0 \in E$ such that $E$ consists of the set of elements
of the form $f\phi_0 + \psi$ with $f \in L^\infty(X)$ and $\psi \in E_0$.
\medskip

Observe that $|\Phi_0(\psi)| \leq |\psi|$ almost everywhere, for any
$\psi \in E_0$. For if not, there would be a positive measure subset $A
\subset X$ such that $|\Phi_0(\psi)| \geq |\psi| + \epsilon$ almost everywhere
on $A$. But then
$$\|\Phi_0(\chi_A\psi)\| = \| |\Phi_0(\chi_A\psi)|\|
= \|\chi_A |\Phi_0(\psi)|\| \geq \|\chi_A\psi\| + \epsilon,$$
contradicting the fact that $\|\Phi_0\| = 1$.
\medskip

The remaining argument is a simple reworking of
the usual proof of the Hahn-Banach theorem. Namely, observe
that if $f_1, f_2 \in L^\infty(X)$, $f_1, f_2 \geq 0$, and $\psi_1, \psi_2 \in E_0$ then
$$\Phi_0(f_2\psi_1 - f_1\psi_2) \leq |f_2\psi_1 - f_1\psi_2|
\leq f_1|f_2\phi_0 + \psi_2| + f_2|f_1\phi_0 + \psi_1|.$$
From this it follows that
$$(-|f_2\phi_0 + \psi_2| - \Phi_0(\psi_2))/f_2
\leq (|f_1\phi_0 + \psi_1| - \Phi_0(\psi_1))/f_1$$
on the common support of $f_1$ and $f_2$. Also both sides lie in the interval
$[-\|\phi_0\|, \|\phi_0\|]$ wherever they are defined.
Taking the supremum of the left side in $L^\infty(X)$ (assigning the default
value $-\|\phi_0\|$ whenever $f_2 = 0$), we obtain $g \in L^\infty(X)$ such that
$$(-|f\phi_0 + \psi| - \Phi_0(\psi))/f \leq g \leq
(|f\phi_0 + \psi| - \Phi_0(\psi))/f$$
holds on the support of $f$, for all $f \in L^\infty(X)$ with $f \geq 0$.
By decomposing an arbitrary $f \in L^\infty(X)$ into its positive and
negative parts, this implies that
$$|fg + \Phi_0(\psi)| \leq |f\phi_0 + \psi|;\eqno{(*)}$$
the inequality is immediate where $f$ is positive, has already been
proven where $f = 0$, and holds in the
negative case by replacing $\psi$ with
$-\psi$. So defining $\Phi(f\phi_0 + \psi) = fg + \Phi_0(\psi)$ for all
$f \in L^\infty(X)$ and $\psi \in E_0$ produces a norm-one module homomorphism
$\Phi$ which extends $\Phi_0$.
\medskip

($\Phi$ is well-defined by the following argument. Suppose $f_1\phi_0 + \psi_1
= f_2\phi_0 + \psi_2$. Then $(f_1 - f_2)\phi_0 + (\psi_1 - \psi_2) = 0$, so
($*$) implies that $|(f_1 - f_2)g + \Phi_0(\psi_1 - \psi_2)| = 0$. Thus
$$(f_1g + \Phi_0(\psi_1)) - (f_2g + \Phi_0(\psi_2))
= (f_1 - f_2)g + \Phi_0(\psi_1 - \psi_2) = 0,$$
which shows that $\Phi$ is well-defined.)\hfill\hal
\bigskip

\noindent {\bf Corollary 6.} {\it If $E$ is an $L^\infty(X)$-normed module
then for any $\phi \in E$ there exists $\Phi \in \E'$ with $\|\Phi\| = 1$ and
$\Phi(\phi) = |\phi|$. The natural map from $E$ into $E''$ is isometric.}
\medskip

\noindent {\it Proof.} The natural map is automatically nonexpansive.
Conversely, given any $\phi \in E$ define $\Phi_0(f\phi) = f|\phi|$;
this is a norm-one module homomorphism defined on $E_0 = \{f\phi:
f \in L^\infty(X)\}$, so it extends to a
norm-one $\Phi \in E'$ by the theorem. Since $\|\Phi(\phi)\| =
\| |\phi| \| = \|\phi\|$, we are done.\hfill\hal
\bigskip
\bigskip

\noindent {\bf B. Characterizations}
\bigskip

First, we characterize abelian W*-modules in terms of duality.
\bigskip

\noindent {\bf Theorem 7.} {\it If $E$ is a Banach module over $L^\infty(X)$
then $E'$ is an abelian W*-module. Conversely, any abelian W*-module over
$L^\infty(X)$ is isometrically and weak*-continuously isomorphic to the dual
of some $L^\infty(X)$-normed module.}
\medskip

\noindent {\it Proof.} To show that $E'$ is an $L^\infty(X)$-normed
module, we check the property given in Theorem 2 (c). Let $p \in L^\infty(X)$
be a projection and let $\Phi \in E'$. Then
$\|p\Phi\| \leq \|p\|\|\Phi\| = \|\Phi\|$, and similarly
$\|(1-p)\Phi\| \leq \|\Phi\|$, so
$$\|\Phi\| \geq {\rm max}(\|p\Phi\|, \|(1-p)\Phi\|).$$
Conversely, given $\epsilon > 0$ find $\phi \in E$ such that $\|\phi\| = 1$
and $\|\Phi(\phi)\| \geq \|\Phi\| - \epsilon$. Then
$$\|\Phi(\phi)\| = {\rm max}(\|p\Phi(\phi)\|, \|(1-p)\Phi(\phi)\|),$$
so we must have
$${\rm max}(\|p\Phi\|, \|(1-p)\Phi\|) \geq
{\rm max}(\|p\Phi(\phi)\|, \|(1-p)\Phi(\phi)\|) = \|\Phi(\phi)\|
\geq \|\Phi\| - \epsilon.$$
This verifies the $L^\infty(X)$-normed property.
\medskip

To see that $E'$ is a dual Banach space, let $Y$ be the set-theoretic
cartesian product of the unit ball of $E$ with the unit ball of $L^1(X)$,
and define $T: E' \to l^\infty(Y)$ by $(T\Phi)(\phi, f) = \int \Phi(\phi)f$.
It is straightforward to check
that $T$ is isometric. If $(\Phi_i)$ is a bounded
universal net in $E'$ then $\Phi(\phi) = \lim \Phi_i(\phi)$ (where the
limit is taken in the weak* topology on $L^\infty(X)$) defines $\Phi \in E'$
such that $T\Phi_i \to T\Phi$ pointwise. This shows that the unit ball of
$T(E')$ is weak*-closed in $l^\infty(Y)$, hence $T(E')$ is weak*-closed by
the Krein-Smulian theorem ([14], Theorem 12.1) and therefore a dual space.
Since $T$ is isometric, $E'$ is also a dual space, and on bounded sets its
weak* topology satisfies $\Phi_i \to \Phi$ if $\Phi_i(\phi) \to \Phi(\phi)$
weak* in $L^\infty(X)$ for all $\phi \in E$.
\medskip

Finally, if $f_i$ is a bounded net in $L^\infty(X)$ which converges weak*
to $f$ then for any $\phi \in E$ and $\Phi \in E'$ we have
$$(f_i\Phi)(\phi) = f_i\cdot \Phi(\phi) \to f\cdot\Phi(\phi) =
(f\Phi)(\phi),$$
so the map $(f, \Phi) \mapsto f\Phi$ is continuous in the first variable.
Whereas if $\Phi_i$ is a bounded net in $E'$ which converges weak* to
$\Phi$, then for any $f \in L^\infty(X)$ and $\phi \in E$ we have
$$(f\Phi_i)(\phi) = f\cdot \Phi_i(\phi) \to f\cdot \Phi(\phi)
= (f\Phi)(\phi),$$
so the map $(f, \Phi) \mapsto f\Phi$ is also continuous in the second variable.
(It is enough to check continuity on bounded nets by the Krein-Smulian
theorem.) So $E'$ is an abelian W*-module, completing the proof of the first
statement.
\medskip

For the converse, let $E$ be an abelian W*-module. By
([21], Theorem 4.1) there is a complex Hilbert space $H$, an isometric
and weak*-continuous map $\imath$ of $E$ onto a weak*-closed subspace of
$B(H)_{sa}$, and an isometric weak*-continuous algebra homomorphism $\pi$
from $L^\infty(X)$ into $B(H)_{sa}$, such that
$$\imath(f)\pi(\phi) = \pi(f\phi) = \pi(\phi)\imath(f)$$
for all $f \in L^\infty(X)$ and $\phi \in E$. Here $B(H)_{sa}$ denotes
the set of bounded self-adjoint operators on $H$.
\medskip

Now let $\phi \in E$ and $\epsilon > 0$; we will find $\Phi \in E'$ which is
weak*-continuous and satisfies $\|\Phi\| = 1$ and $\|\Phi(\phi)\| \geq
\|\phi\| - \epsilon$. By the preceding we may assume $L^\infty(X) \subset
B(H)_{sa}$, $E \subset B(H)_{sa}$,
and $f\phi = \phi f$ for all $f \in L^\infty(X)$ and $\phi \in E$.
\medskip

Since $\phi \in B(H)$ is self-adjoint, we can find a unit vector $v \in H$
such that $|\langle \phi v, v\rangle| \geq \|\phi\| - \epsilon$. Let
$X_0 \subset X$ be the largest subset such that
$\langle f v, v\rangle = 0$ for no $f \in L^\infty(X_0)$ besides $f = 0$,
i.e.\ $X_0$ is the support of the vector state given by $v$.
By a standard argument (e.g.\ see [35]), there is then a unique map
--- a conditional expectation --- $\Phi$ from $E$ into $L^\infty(X_0)
\subset L^\infty(X)$ with the property that
$$\langle \Phi(\psi) f v, v\rangle = \langle \psi f v, v\rangle$$
for all $f \in L^\infty(X)$ and $\psi \in E$, and it is straightforward to
verify that $\Phi \in E'$ has the desired properties.
Thus, letting $E_0$ be the
weak*-continuous part of $E'$, it follows that the natural map from $E$ onto
$E_0'$ is noncontractive. But it is automatically nonexpansive and
weak*-continuous, so $E \cong E_0'$.\hfill\hal
\bigskip

We now deduce a measurable version of Theorem 2 (b).
\bigskip

\noindent {\bf Corollary 8.} {\it There is a finitely decomposable measure
space $Y$, an isometric weak*-continuous algebra homomorphism $\imath$ from
$L^\infty(X)$ into $L^\infty(Y)$, and an isometric weak*-continuous linear map
$\pi$ from $E$ into $L^\infty(Y)$, such that $\imath(f)\pi(\phi) = \pi(f\phi)$
for all $f \in L^\infty(X)$ and $\phi \in E$.}
\medskip

\noindent {\it Proof.} Retain the notation used in the proof of Theorem 7.
For each $\Phi \in E_0$ with $\|\Phi\| = 1$, let
$X_\Phi$ be a copy of $X$; then let $Y = \bigcup X_\Phi$ be their disjoint
union. Let $\imath$ be the diagonal embedding of $L^\infty(X)$ into
$L^\infty(Y)$, and define $\pi: E \to L^\infty(Y)$ by
$\pi(\phi) = \Phi(\phi)$ on $X_\Phi$.
\medskip

It follows from Theorem 7 that $\pi$ is isometric and weak*-continuous,
and the remainder is easy.\hfill\hal
\bigskip

Since $L^\infty(X)$ is semihereditary [28], every finitely-generated
projective module over $L^\infty(X)$ is isomorphic to a direct sum of
ideals (e.g.\ see [11]).
(I am indebted to Ken Goodearl for this argument.) Finitely-generated
abelian W*-modules enjoy a similar characterization.
\bigskip

\noindent {\bf Lemma 9.} {\it For each $n \in \N$ there is a maximal subset
$X_n \subset X$ with the property that $\chi_{X_n}E$ is generated by $n$
elements.}
\medskip

\noindent {\it Proof.} Let $\Gamma$ be the collection of all measurable
subsets $S$ of $X$ with the property that $\chi_SE$ is generated by $n$
elements. (By ``generated'' we mean that $\chi_SE$ is the
smallest abelian W*-module which contains the $n$ generators.) It is clear
that if $S \in \Gamma$ and $T \subset S$ then $T \in \Gamma$. Also,
if $\{S_k\} \subset \Gamma$ and the $S_k$ are disjoint, then let
$\phi_1^k, \ldots, \phi_n^k$ be a generating set for $\chi_{S_k}E$ for
each $k$. Normalizing, we may suppose that $\|\phi_i^k\| \leq 1$ for all
$i$ and $k$. As the $S_k$ are disjoint, we may then
take the weak* sum $\phi_i =
\sum_k \phi_i^k$ ($1 \leq i \leq n$), and the $\phi_i$ then generate
$\chi_S E$ where $S = \bigcup S_k$.
\medskip

Thus, $\Gamma$ is closed under subsets and disjoint unions; it follows that
it contains a maximal element (up to null sets).\hfill\hal
\bigskip

\noindent {\bf Theorem 10.} {\it Let $E$ be a finitely-generated
abelian W*-module over $L^\infty(X)$. Then there is a partition of $X$,
$X = \bigcup_{n = 1}^m X_n$, and for each $x \in X_n$ a norm $\|\cdot\|_x$
on $\R^n$, such that $E$ is isometrically and weak*-continuously isomorphic
to the set of bounded measurable functions $f$ such that
$f|_{X_n}$ takes $X_n$ into $\R^n$ for
$1 \leq n \leq m$, with norm given by
$\|f\| = {\rm ess}\sup \|f(x)\|_x$.}
\medskip

\noindent {\it Proof.} The differences of the sets described in Lemma 9
provide the desired partition of $X$. By Theorem 2 (c), it suffices to
restrict attention to one block of the partition; thus we may suppose
that $E$ is generated by $n$ elements $\phi_1, \ldots, \phi_n$ and for
any positive measure subset $A \subset X$, $\chi_A E$ is not generated by
fewer than $n$ elements.
\medskip

For $x \in X$ and $a = (a_1, \ldots, a_n) \in \R^n$, define
$$\|a\|_x = |a_1\phi_1 + \cdots + a_n\phi_n|(x).$$
This is a seminorm for almost every $x$ by ([51], Proposition 2), and
the module of bounded measurable sections of the trivial bundle with
fiber $\R^n$, equipped with this family of seminorms, is
isometric to a weak*-dense subset of $E$ via the identification of
$f: X \to \R^n$ with $\sum (\pi_i\circ f)\phi_i$, where $\pi_i$ is the
$i$th coordinate projection on $\R^n$.
Furthermore, the weak*-topology of $E$ agrees with the weak*-topology of
$L^\infty(X)$ on the $L^\infty(X)$-span of each $\phi_i$, hence the
identification just described is a weak*-homeomorphism, hence the module
of sections is isometric to all of $E$ since its unit ball is already
weak* compact.
\medskip

Finally, suppose $\|\cdot\|_x$ fails to be a norm on a set $A$ of positive
measure. For each $x \in A$ let $B_x = \{a \in \R^n: \|a\|_x = 0\}$; then
$B = \bigcup B_x$ is an (F) Banach bundle over $L^\infty(A)$ and so there
exists a bounded, nonzero, measurable section $f: A \to \R^n$ with
$f(x) \in B_x$ for all $x$ by ([23], Appendix). Let $\phi =
\sum (\pi_i \circ f)\phi_i$, setting $\phi = 0$ off of $A$;
without loss of generality we may assume $|\pi_1 \circ f| \geq \epsilon > 0$
on a positive measure subset $A_0$ of $A$.
Then $\phi_1$ is expressible as a linear combination of $\phi$
together with $\phi_2, \ldots, \phi_n$. But $|\phi|(x) = 0$ for almost
every $x \in A$, hence $\phi = 0$, so that $\phi_2, \ldots, \phi_n$
generate $\chi_{A_0}E$. This contradicts the reduction made in the first
paragraph of the proof, so we conclude that
$\|\cdot\|_x$ is a norm for almost every $x \in X$.\hfill\hal
\bigskip

Theorem 10 can be viewed as saying that any finitely-generated abelian
W*-module is isometrically isomorphic to the module of bounded measurable
sections of some bundle of finite-dimensional Banach spaces. An analogous
result holds for non-finitely-generated Hilbert modules ([49], Theorem 3.12;
see also [65]); in this case the fibers are Hilbert spaces.
\medskip

Morally, it should be true without the finitely-generated assumption
that every abelian W*-module is isometric to
the module of bounded measurable sections of some bundle of dual Banach
spaces. However, measure-theoretic complications make it difficult to formulate
a satisfying general result of this type.
\medskip

We require one final construction: the tensor product of abelian
W*-modules. For our purposes the appropriate definition is the following.
\bigskip

\noindent {\bf Definition 11.} Let $E = E_0'$ and $F = F_0'$ be abelian
W*-modules over $L^\infty(X)$. We define $E\otimes F =
E\otimes_{L^\infty(X)} F$ to be the set of bounded module maps
from $E_0$ into $F$, or equivalently the set of bounded module maps from
$F_0$ into $E$. It is straightforward to check that this is again an
abelian W*-module.
\medskip

(A bounded module map $T: E_0 \to F$ gives rise to a map from $F_0 \subset
F'$ to $E$ by taking adjoints, and vice versa.)\hfill\hal
\bigskip
\bigskip

\noindent {\underbar{\bf 3. Measurable metrics.}}
\bigskip

\noindent {\bf A. Definitions}
\bigskip

In some ways the pointwise aspect of metrics does not interact well with
non-atomic measures. A helpful alternative, which involves only distances
between positive measure sets, was formulated in [68] and also used in [69].
It is the following.
\bigskip

\noindent {\bf Definition 12.} Let $M = (M, \mu)$ be a $\sigma$-finite
measure space and let $\Omega$ be the collection of all positive measure
subsets of $M$, modulo null sets. A {\it measurable pseudometric} is a
map $\rho: \Omega^2 \to [0, \infty]$ such that
$$\eqalign{\rho(A, A) &= 0\cr
\rho(A, B) &= \rho(B, A)\cr
\rho(\bigcup_{n=1}^\infty A_n, B) &= \inf_n \rho(A_n, B)\cr
\rho(A, C) &\leq \sup_{B' \subset B} (\rho(A, B') + \rho(B', C))\cr}$$
for all $A, B, C, A_n \in \Omega$.\hfill\hal
\bigskip

\noindent {\bf Definition 13.} Let $M = (M, \mu, \rho)$ be a $\sigma$-finite
measure space with measurable pseudometric. The {\it essential range} of
$f \in L^\infty(M)$ is the set of $a \in \R$ such that $f^{-1}(U)$ has
positive measure for every neighborhood $U$ of $a$ (equivalently, it is the
spectrum of $f$ in the real Banach algebra $L^\infty(M)$); and we let
$\rho_f(A, B)$ denote the distance, in $\R$, between the essential ranges of
$f|_A$ and $f|_B$.
\medskip

The {\it Lipschitz number} $L(f)$ of $f$ is
$$L(f) = \sup \{\rho_f(A, B)/\rho(A, B): A, B \in \Omega\hbox{ and }\rho(A, B)> 0\}$$
and $\Li(M)$ is the set of all $f \in L^\infty(M)$ for which $L(f)$ is
finite. With the norm $\|f\|_L = {\rm max}(L(f), \|f\|_\infty)$, $\Li(M)$
is easily seen to be a Banach space as well as a ring.\hfill\hal
\bigskip

We say that $\rho$ is a {\it measurable metric} if $\Li(M)$ is weak*-dense
in $L^\infty(M)$. We also have the following equivalent condition, which,
in the case of an atomic measure, is equivalent to the condition
$\rho(x,y) = 0 \Rightarrow x = y$.
\bigskip

\noindent {\bf Proposition 14.} {\it A measurable pseudometric $\rho$ is a
measurable metric if and only if the underlying measurable $\sigma$-algebra is
generated (up to null sets) by the sets $A \in \Omega$ with the property
that
$$A \cap B = \emptyset\qquad \Rightarrow \qquad
\hbox{there exists }B' \subset B\hbox{ such that }\rho(A, B')>0$$
for all $B \in \Omega$.}
\medskip

\noindent {\it Proof.} Call a set $A$ {\it full} if it satisfies the
displayed condition. If the full sets do not generate the measurable
$\sigma$-algebra then $L^\infty(M)$ strictly contains the weak*-closed
algebra generated by the characteristic functions $\chi_A$ for
$A \in \Omega$ full. However, we claim that $\Li(M)$ is contained in this
weak*-closed algebra, so that failure of the displayed condition implies
that $\rho$ is not a measurable metric. To verify the claim observe that if
$f \in \Li(M)$ and $a, b \in \R$, $a \leq b$, then $A = f^{-1}([a, b])$ is
full. For if $B \in \Omega$ is disjoint from $A$ then the essential range of
$f|_B$ is not contained in $[a, b]$, hence
$$B' = B \cap f^{-1}((-\infty, a - \epsilon] \cup [b+\epsilon, \infty))$$
has positive measure for some $\epsilon > 0$, and $\rho_f(A, B') \geq
\epsilon$ hence $\rho(A, B') \geq \epsilon/L(f) > 0$.
From this it follows that $f$ can be approximated in sup norm
by simple functions which are measurable with respect to the $\sigma$-algebra
generated by the full sets. This completes the proof of the forward
direction.
\medskip

Conversely, suppose the full sets do generate the measurable $\sigma$-algebra.
For any full set $A$ define $f_A \in \Li(M)$ by
$$f_A = \sup_{B \in \Omega} (\rho(A, B) \wedge 1)\cdot \chi_B$$
where $\rho(A, B) \wedge 1 = {\rm min}(\rho(A, B), 1)$ and the supremum
is taken in $L^\infty(M)$. To see that $f \in \Li(M)$
let $B, C \in \Omega$ and $\epsilon > 0$. By ([68], Lemma 5) there exist
$B' \subset B$ and $C' \subset C$ such that $\rho(B'', C'') \leq \rho(B, C) +
\epsilon$ for all $B'' \subset B'$ and $C'' \subset C'$. Applying
([68], Lemma 5) again we may assume that $\rho(A, B'') \leq \rho(A, B') + \epsilon$
and $\rho(A, C'') \leq \rho(A, C') + \epsilon$ for all $B'' \subset B$ and
$C'' \subset C$; this implies that the essential range of $f_A|_B$
intersects the interval $[\rho(A, B'), \rho(A, B') + \epsilon]$ and the
essential range of $f_A|_C$ intersects the interval
$[\rho(A, C'), \rho(A, C') + \epsilon]$. Without loss of generality suppose
$\rho(A, B') \leq \rho(A, C')$. Then
$$\rho_{f_A}(B, C) \leq \rho(A, C') - \rho(A, B') + \epsilon,$$
and finding $B'' \subset B'$ such that $\rho(A, C') \leq \rho(A, B'') +
\rho(B'', C') + \epsilon$ we then have
$$\eqalign{\rho_{f_A}(B, C)
&\leq (\rho(A, B'') + \rho(B'', C') + \epsilon) - \rho(A, B') + \epsilon\cr
&\leq \rho(B, C) + 4\epsilon.\cr}$$
We conclude that $L(f_A) \leq 1$.
\medskip

Now $f^{-1}(\{0\}) = A$ up to a null set, so the weak*-closed subalgebra of
$L^\infty(M)$ generated by $\Li(M)$ contains $\chi_A$. We have therefore
shown that the characteristic function of any full
set belongs to this algebra, hence the algebra equals $L^\infty(M)$ and so
$\rho$ is a measurable metric.\hfill\hal
\bigskip

We record two basic facts. The first follows from ([69], Theorem 9) and
the second is ([67], Theorem B).
(The latter was proven in [67] only for pointwise metrics, but the proof in
the measurable case is essentially identical.)
\bigskip

\noindent {\bf Proposition 15.} {\it $\Li(M)$ is a dual space, and on its
unit ball the weak*-topology agrees with the restriction of the
weak*-topology on $L^\infty(M)$.}\hfill\hal
\bigskip

\noindent {\bf Theorem 16.} {\it Let $M$ be a measurable metric space and
let $\A$ be a weak*-closed subalgebra of $\Li(M)$. Suppose that
there exists $k \geq 1$ such that for every $A, B \subset M$ we have
$$\rho_f(A, B) = \rho(A, B)$$
for some $f \in \A$ with $L(f) \leq k$. Then $\A = \Li(M)$.}\hfill\hal
\bigskip

\noindent {\bf B. Derivations}
\bigskip

Measurable metrics are closely connected to a natural class of derivations.
These derivations are described in the following definitions.
\bigskip

\noindent {\bf Definition 17.} Let $E$ be a bimodule over $L^\infty(X)$
(with possibly different left and right actions) which is also a dual
Banach space. It is an {\it abelian W*-bimodule} if there is a finitely
decomposable measure space $Y$, an isometric weak*-continuous linear map
$\pi$ from $E$ into $L^\infty(Y)$, and isometric weak*-continuous algebra
homomorphisms $\imath_l$ and $\imath_r$ from $L^\infty(X)$ into
$L^\infty(Y)$ such that
$$\pi(f\phi g) = \imath_l(f)\pi(\phi)\imath_r(g)$$
for all $f, g \in L^\infty(X)$ and $\phi \in E$.\hfill\hal
\bigskip

\noindent {\bf Definition 18.} Let $E$ be an abelian W*-bimodule over
$L^\infty(X)$. An (unbounded)
{\it W*-derivation} from $L^\infty(X)$ into $E$ is then
a linear map $\delta$ from a weak*-dense, unital subalgebra of $L^\infty(X)$
into $E$ with the property that $\delta(fg) = f\delta(g) + \delta(f)g$
and whose graph is a weak*-closed subspace of $L^\infty(X) \oplus E$.\hfill\hal
\bigskip

The following theorem is a slight reformulation of
the main result of [69]. We will use it in section 6.
\bigskip

\noindent {\bf Theorem 19.} {\it Let $M = (X, \rho)$ be a measurable metric
space. Then there is a W*-derivation
$\delta$ from $L^\infty(X)$ into an abelian W*-bimodule whose domain
equals $\Li(M)$ and such that $L(f) = \|\delta(f)\|$ for all $f \in
\Li(M)$.
\medskip

Conversely, let $\delta$ be any W*-derivation from $L^\infty(X)$ into an
abelian W*-bimodule. Then there is a measurable metric $\rho$ on $X$ such that
the domain of $\delta$ equals $\Li(M)$ and $L(f) = \|\delta(f)\|$ for all
$f \in \Li(M)$, where $M = (X, \rho)$.}\hfill\hal
\bigskip

\noindent {\bf C. Metric realization}
\bigskip

Any genuine metric $\rho$ on $M$ gives rise to a measurable metric, by setting
$$\rho_0(A,B) = \inf_{x \in A, y \in B} \rho(x,y)$$
and then letting $\rho(A, B)$ be the supremum of $\rho_0(A', B')$ as $A'$ and
$B'$ range over all measurable sets which differ from $A$ and $B$ by null
sets. If $\mu$ is atomic, it is not hard to see that every measurable metric
on $M$ comes from a unique pointwise metric in this manner.
\medskip

But in general not every measurable metric on $M$ arises in this way, and
in the second part of Theorem 19 the use of measurable metrics is
necessary. However, if one is willing to modify the set $X$ one can always
get a genuine underlying metric. This is in keeping with the algebraic
point of view which regards the algebra $L^\infty(X)$ as primary and the
measure space $X$ as non-canonical and secondary. We now show how an
arbitrary measurable metric can be reduced to an ordinary metric; this
incidentally sharpens the results of [68] and [69].
\bigskip

\noindent {\bf Theorem 20.} {\it Let $M = (M, \mu, \rho_M)$ be a
$\sigma$-finite
measurable metric space. Then there is a $\sigma$-finite measure space $N =
(N, \nu)$ with a complete pointwise metric $\rho_N$ and an isometric
isomorphism of $L^\infty(M)$ onto $L^\infty(N)$ which carries $\Li(M)$
onto $\Li(N)$ ( = the bounded measurable Lipschitz functions on $N$)
in a manner which preserves Lipschitz number.}
\medskip

\noindent {\it Proof.} For the duration of the proof we switch to complex
scalars. Let $\A$ be the C*-subalgebra of $L^\infty(M)$ generated by
$\Li(M)$, and let $N \subset \A'$
be the spectrum of $\A$, so that we have a canonical
identification of $\A$ with $C(N)$. Choose a nowhere-zero function
$f \in L^1(M)$ and let $\nu$ be the Borel measure on $N$ which represents the
linear functional on $\A$ given by integrating against $f$. It is standard
that $L^\infty(M)$ is then canonically isometrically isomorphic to
$L^\infty(N)$. For $f \in L^\infty(M)$ let $\tilde f$ denote the
corresponding function in $L^\infty(N)$.
\medskip

For $\phi, \psi \in N$ define
$$\rho_N(\phi, \psi) = \sup\{|(\phi - \psi)(f)|: f \in \Li(M), L(f) \leq 1\}.$$
It is straightforward to check that this is a metric on $N$ (possibly
with infinite distances). To verify completeness, let $(\phi_n)$ be a
Cauchy sequence in $N$. Then $(\phi_n)$ is also Cauchy as a sequence in
$\Li(M)'$, so it converges to some $\phi \in \Li(M)'$. As each $\phi_n$
is a complex homomorphism so is $\phi$, and this implies that
$|\phi(f)| \leq \|f\|_\infty$ for all $f \in \Li(M)$; thus $\phi$ extends
by continuity to $\A$, so that $\phi \in N$. Since $\phi_n \to \phi$
weak* in $\Li(M)'$, it is the case that $\rho_N(\phi_m, \phi_n) \leq \epsilon$
for all $m \geq n$ implies $|(\phi - \phi_n)(f)| \leq \epsilon$ for all
$f \in \Li(M)$ with $L(f) \leq 1$; so $\phi_n \to \phi$ in $N$.
\medskip

We must now show that $\Li(M)$ is isometrically identified with $\Li(N)$.
First, let $f \in \Li(M)$ and suppose $L(f) \leq 1$. Then
$$|\tilde f(\phi) - \tilde f(\psi)| = |(\phi - \psi)(f)| \leq \rho_N(\phi, \psi)$$
for any $\phi, \psi \in N$, so that $L(\tilde f) \leq 1$. This shows that
$\Li(M) \subset \Li(N)$ and $L(\tilde f) \leq L(f)$ for any $f \in \Li(M)$.
\medskip

Conversely, let $\tilde f \in \Li(N)$ and suppose $L(\tilde f) \leq 1$.
Let $A, B \subset M$ be positive measure sets with $\rho_M(A, B) > 0$; we
must show that the distance between the essential ranges of $f|_A$ and
$f|_B$ is at most $\rho_M(A, B)$. Let $\epsilon > 0$; by ([68], Lemma 5)
there exist $A' \subset A$ and $B' \subset B$ such that
$\rho_M(A'', B'') \leq \rho_M(A, B) + \epsilon$ for every $A'' \subset A'$ and
$B'' \subset B'$. Let
$\phi$ be a complex homomorphism on $\A$ which factors through restriction to
$A'$, and let $\psi$ be a complex homomorphism on $\A$ which factors
through restriction to $B'$. We claim that $\rho_N(\phi, \psi) \leq
\rho_M(A, B) + \epsilon$.
\medskip

To verify the claim, suppose it fails; then there exists $g \in \Li(M)$
with $L(g) \leq 1$ such that $|(\phi - \psi)(g)| > \rho_M(A, B) + \epsilon$.
Let $U$ and $V$ be open neighborhoods of $\phi(g)$ and $\psi(g)$ whose
distance also exceeds $\rho_M(A, B) + \epsilon$. Then $\phi(g)$ and $\psi(g)$
belong to the essential ranges of $g|_{A'}$ and $g|_{B'}$,
so $g^{-1}(U) \cap A'$
and $g^{-1}(V) \cap B'$ are positive measure sets whose distance in $M$ is
greater than $\rho_M(A, B) + \epsilon$, since $L(g) \leq 1$. But this
contradicts the choice of $A'$ and $B'$, so the claim is proven.
\medskip

Finally, $\phi(f)$ and $\psi(f)$ belong to the essential ranges of $f|_A$
and $f|_B$, so
$$\rho_f(A, B) \leq |(\phi - \psi)(f)| \leq \rho_N(\phi, \psi) \leq
\rho_M(A, B) + \epsilon.$$
Taking $\epsilon$ to zero, we conclude that $\rho_f(A,B) \leq
\rho_M(A, B)$. This shows that $\Li(N) \subset \Li(M)$ and
$L(f) \leq L(\tilde f)$, completing the proof that $\Li(M)$ is isometrically
identified with $\Li(N)$.\hfill\hal
\bigskip
\bigskip

\noindent {\underbar{\bf 4. 1-forms and the exterior derivative.}}
\bigskip

\noindent {\bf A. The construction}
\bigskip

\noindent {\bf Definition 21.} Let $M$ be a measurable metric space and
let $E$ be an abelian W*-module over $L^\infty(M)$. A {\it metric derivation}
$\delta: \Li(M) \to E$ is a bounded weak*-continuous linear map which
satisfies the derivation identity $\delta(fg) = f\delta(g) + \delta(f)g$
for all $f,g \in \Li(M)$.
\medskip

The {\it module of measurable vector fields $\X(M)$}
is the set of all metric derivations
$\delta: \Li(M) \to L^\infty(M)$. The module action is given by
$(f\cdot\delta)(g) = f\delta(g)$ for $f \in L^\infty(M)$, $\delta \in \X(M)$,
and $g \in \Li(M)$. Using Theorem 2 (c),
it is straightforward to check that $\X(M)$ is an $L^\infty(M)$-normed module.
\medskip

The {\it module of measurable 1-forms $\Omega(M)$} is the dual module of
$\X(M)$, that is $\Omega(M) = \X(M)'$. By Theorem 7 this is an abelian
W*-module over $L^\infty(M)$.\hfill\hal
\bigskip

In general $\X(M)$ is not an abelian W*-module over $L^\infty(M)$; see
Proposition 45 (a). This is true in the finitely-generated and Hilbert
module cases, however
(Corollary 24). Also, the notion of metric derivation is strictly
weaker than the notion of W*-derivation given in Definition 18. In particular,
there may exist elements of $\X(M)$
which are not W*-derivations; see Proposition 45 (b).
\medskip

By Theorem 10, if $\Omega(M)$ is finitely-generated then it can be realized
as the module of bounded measurable sections of some bundle of
finite-dimensional Banach spaces. The latter plays the role of the cotangent
bundle. Again, something like this possibly involving infinite-dimensional
dual Banach spaces should be true even when $\Omega(M)$ is not
finitely-generated. These would act as cotangent spaces and their preduals
as tangent spaces. In the finitely-generated case we assuredly have a
natural tangent bundle according to Corollary 24 and Theorem 10.
\medskip

Additionally, in the Hilbert module case we have tangent and
cotangent bundles, even if $M$ is infinite-dimensional. Here they are
bundles of Hilbert spaces, and the tangent and cotangent spaces at each
point are naturally identified with each other (Corollary 24; see the
comment following Theorem 10).
\medskip

We should remark that Gromov has given a definition of the tangent space
or ``asymptotic cone''
at a point of any metric space [30]. It appears to bear little
relationship to our definition. For instance, the Gromov tangent space
at a boundary point of a manifold with boundary will be a half-space,
whereas our tangent spaces are always Banach spaces. On the other hand,
our tangent spaces will only be well-defined almost everywhere.
\medskip

By taking tensor products over $L^\infty(M)$, one can define
the module of bounded measurable tensor fields of arbitrary type $(r,s)$.
Together with the following definition this raises the question of whether
higher-order exterior derivatives must exist, i.e.\ whether $d$ can in general
be extended to the whole exterior algebra. I suspect the answer is no
even in the finite-dimensional case, but
this has been done for Lipschitz manifolds [27].
\bigskip

\noindent {\bf Definition 22.} The {\it exterior derivative} on $M$ is
the map $d: \Li(M) \to \Omega(M)$ given by $(df)(\phi) = \phi(f)$ for
$f \in \Li(M)$ and $\phi \in \X(M)$.
\medskip

More generally, if $E$ is any submodule of $\X(M)$ there is a natural map
$d_E: \Li(M) \to E'$ given by the same formula, $(d_Ef)(\phi) = \phi(f)$.
Letting $T: \Omega(M) \to E'$ be the natural projection, we have
$d_E = T\circ d$.\hfill\hal
\bigskip

\noindent {\bf Theorem 23.} {\it The exterior derivative $d$ is a metric
derivation. It is universal in the sense that if $\delta: \Li(M) \to E$ is
any metric derivation into an abelian W*-module
then there is a bounded weak*-continuous
$L^\infty(M)$-module map $T: \Omega(M) \to E$ such that $\delta = T\circ d$.
Furthermore, $\|T\| = \|\delta\|$.}
\medskip

\noindent {\it Proof.} It is clear that $d$ is linear. If $f, g \in \Li(M)$
and $\phi \in \X(M)$ then
$$d(fg)(\phi) = \phi(fg) = f\phi(g) + \phi(f)g
= f(dg)(\phi) + (df)(\phi)g,$$
so $d$ is a derivation. And if $(f_i)$ is a bounded net in $\Li(M)$ and
$f_i \to f$ weak* then
$$(df_i)(\phi) = \phi(f_i) \to \phi(f) = (df)(\phi)$$
weak* in $L^\infty(M)$ for any $\phi \in \X(M)$, hence $df_i \to df$ weak*
in $\Omega(M)$. This shows that $d$ is a metric derivation.
\medskip

Let $\delta: \Li(M) \to E = E_0'$ be any metric derivation. (We may assume that
$E$ is a dual module by Theorem 7.) Define $T_0: E_0 \to \X(M)$ by
$(T_0 \phi)(f) = (\delta f)(\phi)$ for $\phi \in E_0$ and $f \in \Li(M)$. Then
$T_0$ is a bounded $L^\infty(M)$-module map and it has a bounded,
weak*-continuous adjoint $T: \Omega(M) \to E$.
For any $f \in \Li(M)$ and $\phi \in E$ we then have
$$T(df)(\phi) = (df)(T_0\phi) = (T_0\phi)(f) = (\delta f)(\phi),$$
so $\delta = T \circ d$. Also $\|T\| = \|T_0\| = \|\delta\|$.\hfill\hal
\bigskip

\noindent {\bf Corollary 24.} {\it If $\Omega(M)$ is reflexive then $\X(M)$
is an abelian W*-module. In particular, this holds if $\X(M)$ is
finitely-generated or if $\X(M)$ satisfies the parallelogram law
$$|\phi + \psi|^2 + |\phi - \psi|^2 = 2|\phi|^2 + 2|\psi|^2$$
(almost everywhere, for all $\phi, \psi \in \X(M)$). In the latter case $\X(M)$
and $\Omega(M)$ are canonically isomorphic self-dual Hilbert modules.}
\medskip

\noindent {\it Proof.} Suppose $\Omega(M)$ is reflexive. By this
we mean that the natural map from $\Omega(M)$ into $\Omega(M)''$ ---
which is isometric by Corollary 6 --- is onto.
This implies that every element of $\Omega(M)'$ is weak*-continuous.
So for any $\Phi \in \Omega(M)'$ the map $f \mapsto \Phi(df)$ is a
metric derivation from $\Li(M)$ into $L^\infty(M)$, i.e.\ it is an element of
$\X(M)$. This shows that the natural map from $\X(M)$ into $\Omega(M)'$ ---
which, again, is isometric by Corollary 6 --- is onto, hence $\X(M)$ is
an abelian W*-module by Theorem 7.
\medskip

Suppose $\X(M)$ is finitely-generated. Then $E = \X(M)''$ is a
finitely-generated abelian W*-module, and Theorem 10 implies that
$E'$ is finitely-generated as well (namely, it is the module of bounded
measurable sections of the same vector bundle, equipped with the fiberwise
dual norms). But $\Omega(M)$ is a quotient of $E' = \Omega(M)''$ by
Theorem 5, so $\Omega(M)$ is also finitely-generated.
By Theorem 10 we deduce that it is reflexive.
\medskip

Now suppose $\X(M)$ satisfies the paralellogram law. Then $\X(M)$ is a
Hilbert module by ([51], Lemma 13) and so $\Omega(M)$ is self-dual and hence
reflexive by ([49], Theorem 3.2). From the first part of the corollary it
now follows that $\X(M) \cong \Omega(M)' \cong \Omega(M)$.\hfill\hal
\bigskip

In the finitely-generated case we have a simple formula
for the module of vector fields on a product space. I suspect it is
false in general, but may be true if $\X(M)$ and $\X(N)$ are Hilbert
modules.
\bigskip

\noindent {\bf Theorem 25.} {\it Let $M$ and $N$ be measurable metric spaces
and suppose that $\X(M)$ and $\X(N)$ are finitely-generated. Then
$$\X(M\times N) \cong \big(\X(M)\otimes L^\infty(M\times N)\big)
\oplus \big(\X(N)\otimes L^\infty(M\times N)\big),$$
where the first tensor product is taken over $L^\infty(M)$ and the second is
taken over $L^\infty(N)$, and $\cong$ denotes
isomorphism of $L^\infty(M\times N)$-modules.}
\medskip

\noindent {\it Proof.} Recall our version of tensor products of modules
given in Definition 11. Here we are viewing $L^\infty(M\times N)$ as
an abelian W*-module over either $L^\infty(M)$ or $L^\infty(N)$.
\medskip

We need not specify the product metric on $M\times N$
exactly, since all natural choices
are bi-Lipschitz equivalent. This ambiguity corresponds to a
choice of the direct sum norm in the right side of the isomorphism.
We do require that $\rho(A\times B, A'\times B) = \rho(A, A')$ and
$\rho(A\times B, A\times B') = \rho(B, B')$ for all $A, A' \subset M$ and
$B, B' \subset N$. The easiest way to see that measurable metrics on
$M\times N$ satisfying these conditions exist is via Theorem 20.
\medskip

We first define a map $S$ from $\X(M\times N)$ into the right side. To do
this it suffices to separately define maps $S_M: \X(M\times N) \to
\X(M)\otimes L^\infty(M\times N)$ and $S_N: \X(M\times N) \to
X(N)\otimes L^\infty(M\times N)$. Since $\Omega(M)$ is
reflexive, Corollary 24 implies that $\X(M) \cong \Omega(M)'$, and so
$\X(M)\otimes L^\infty(M\times N)$ is identified with the set of
bounded $L^\infty(M)$-module maps from $\Omega(M)$ into $L^\infty(M\times N)$.
Thus, to define $S_M$ we must say how to produce such a map from an
arbitrary derivation $\delta \in \X(M\times N)$. This is done by observing
that $\Li(M)$ naturally imbeds in $\Li(M\times N)$, hence $\delta$
restricts to a metric derivation from $\Li(M)$ into $L^\infty(M\times N)$,
and the required module map is then given by the universality statement
in Theorem 23. $S_N$ is defined similarly.
It is clear that $\|S_M\|, \|S_N\| \leq 1$.
\medskip

Next, we show that $S(\X(M\times N))$ contains $\X(M)\otimes 1_{M\times N}$
and $\X(N)\otimes 1_{M\times N}$. Identifying these sets with $\X(M)$ and
$\X(N)$ in the obvious way, we will define maps
$T_M: \X(M) \to \X(M\times N)$ and $T_N: \X(N) \to \X(M\times N)$ such
that $S\circ T_M$ is the identity on $\X(M)$ and $S\circ T_N$ is
the identity on $\X(N)$.
\medskip

Let $\delta \in \X(M)$ and let $f \in \Li(M\times N)$; to define
$T_M(\delta)$ we must produce an element of $L^\infty(M\times N)$,
or equivalently a bounded linear map from $L^1(N)$ into $L^\infty(M)$.
Letting $g \in L^1(N)$, the desired element of $L^\infty(M)$ is then
$\delta(f_g)$ where
$$f_g = \int f(x,y) g(y)dy$$
is a bounded Lipschitz function with $L(f_g) \leq L(f)\|g\|_1$. It is
now a matter of unwrapping definitions to verify that
$S_M(T_M(\delta)) = \delta$ and $S_N(T_M(\delta)) = 0$. A similar argument
proves the corresponding result for $T_N$.
\medskip

Now $\X(M)\otimes 1_{M\times N}$ and $\X(N)\otimes 1_{M\otimes N}$ together
algebraically generate the target space as an $L^\infty(M\times N)$-module.
This can be seen by using the
explicit structure of $\X(M)$ and $\X(N)$ given by Corollary 24 and
Theorem 10, and it implies from the above that $S$ is onto. $S$ is also
1-1 since the algebraic tensor product of $\Li(M)$ and $\Li(N)$ is
weak*-dense in $\Li(M\times N)$ by Theorem 16.
So the open mapping theorem implies that $S$ is an isomorphism.\hfill\hal
\bigskip

We conclude this section with a technical criterion which is helpful in
determining $\X(M)$ in some examples.
\bigskip

\noindent {\bf Theorem 26.} {\it Let $M$ be a measurable metric space.
Let $E$ be a submodule of $\X(M)$ which is reflexive as an
$L^\infty(M)$-module and suppose $|df| = |d_Ef|$ for all $f$ in a
weak*-dense subspace of $\Li(M)$. Then $E = \X(M)$.}
\medskip

\noindent {\it Proof.} Let $S \subset \Li(M)$ be a weak*-dense subspace
such that $|df| = |d_E f|$ for all $f \in S$. Let $T: \Omega(M) \to E'$
be the restriction map, so that $d_E = T \circ d$ and $T$ is nonexpansive.
\medskip

Moreover, we have $|T(df)| = |d_E f| = |df|$ for all $f \in S$. We now
claim that $T$ remains isometric on the $L^\infty(M)$-span of these
elements $df$.
To see this let $\sum_1^n f_idg_i$ be a finite linear combination with
$f_i \in L^\infty(M)$ and $g_i \in S$, and
let $\epsilon > 0$. Without loss of generality suppose $L(g_i) \leq 1$ for
all $i$. Let $A$ be a positive measure set on which
$$|\sum f_idg_i| \geq \|\sum f_idg_i\| - \epsilon.$$
By shrinking $A$, we may assume that each $f_i$ varies by at most
$\epsilon/n$ on $A$. Choose $a_i \in \R$ such that
$|f_i(x) - a_i| \leq \epsilon/n$ for almost every $x \in A$.
\medskip

Now
$$|(\sum f_idg_i) - (\sum a_idg_i)| \leq
\sum |(f_i - a_i)dg_i| \leq \sum |f_i - a_i|L(g_i) \leq \epsilon$$
on $A$, so
$$|d(\sum a_ig_i)| = |\sum a_i dg_i|
\geq |\sum f_i dg_i| - \epsilon \geq \|\sum f_idg_i\| - 2\epsilon$$
on $A$. Since $|T(df)| = |df|$ for all $f \in S$, we then have
$$|\sum a_i d_E g_i| = |d_E(\sum a_ig_i)|
= |d(\sum a_ig_i)| \geq \|\sum f_idg_i\| - 2\epsilon$$
on $A$. But finally $|(\sum f_id_E g_i) - (\sum a_i d_E g_i)|
\leq \epsilon$ on $A$ by applying $T$ to the inequality
$|(\sum f_i dg_i) - (\sum a_i dg_i)| \leq \epsilon$, and so we conclude that
$$\|\sum f_id_E g_i\| \geq \|\sum f_idg_i\| - 3\epsilon.$$
Taking $\epsilon$ to zero completes the proof of the claim.
\medskip

Let $\phi \in \X(M)$. Let $E_0$ be the set of elements in $E'$
of the form $\sum f_id_E g_i$ considered above,
and let $T^{-1}$ denote the isometric
embedding of $E_0$ into $\Omega(M)$ given by $T^{-1}(\sum f_id_E g_i) =
\sum f_i dg_i$. Then the map $\Phi \mapsto T^{-1}(\Phi)(\phi)$ is a
bounded module homomorphism from $E_0$ into $L^\infty(M)$, and by Theorem 5
this extends to an element of $E'' = E$. Thus there exists $\phi_0 \in E$
such that
$$\phi(f) = (df)(\phi) = T^{-1}(d_E f)(\phi) = (d_E f)(\phi_0)
= \phi_0(f)$$
for all $f \in S$, hence $\phi = \phi_0$. So $\X(M) = E$.\hfill\hal
\bigskip

It is natural to conjecture that the hypothesis
$\|d_E f\| = \|df\|$ for all $f \in \Li(M)$ should imply the condition
$|d_E f| = |df|$ needed in the preceding theorem. However, even if $M$
is differentiable in the sense of Definition 30, this implication is
false; see Theorem 55.
\bigskip
\bigskip

\noindent {\bf B. Locality}
\bigskip

Since we are treating abelian W*-modules as bimodules with identical left and
right actions, our metric derivations have a local character that contrasts
with the bimodule derivations considered in [69]. This is seen in the
following results.
\bigskip

\noindent {\bf Lemma 27.} {\it Let $\delta: \Li(M) \to E$ be a metric
derivation and let $A \subset M$. If $f, g \in \Li(M)$ satisfy
$f|_A = g|_A$ almost everywhere then
$(\delta f)|_A = (\delta g)|_A$ almost everywhere.}
\medskip

\noindent {\it Proof.} The statement that $(\delta f)|_A = (\delta g)|_A$ is
to be interpreted as meaning that $\chi_A\cdot \delta f = \chi_A \cdot
\delta g$. By considering $f - g$, it will suffice to show that
$f|_A = 0$ implies $(\delta f)|_A = 0$.
Let $I = \{f \in \Li(M): f|_A = 0\}$. This
is a weak*-closed ideal of $\Li(M)$, so by ([69], Theorems 3 and 9) $I^2$ is
weak*-dense in $I$. Thus, for any $f \in I$ we can find a pair of nets
$(f_i), (g_i) \subset I$ such that $f_ig_i \to f$ weak*. But then
$$f_i\delta (g_i) + \delta (f_i)g_i = \delta(f_ig_i) \to \delta(f)$$
weak* in $E$, and since $f_i, g_i \in I$ we have
$\delta (f_ig_i)|_A = 0$ for all $i$, hence $(\delta f)|_A = 0$.\hfill\hal
\bigskip

\noindent {\bf Lemma 28.} {\it Let $M$ be a measurable metric space, let
$A \subset M$, and let $f \in \Li(A)$. Then there exists $g \in \Li(M)$
such that $g|_A = f$ and $L(g) = L(f)$.}
\medskip

\noindent {\it Proof.} Without loss of generality suppose $\|f\|_L = 1$.
For any positive measure set $B \subset A$, define $f_B \in \Li(M)$ as
in the proof of Proposition 14 by
$$f_B = \sup_C (\rho(B, C) \wedge 2)\cdot \chi_C.$$
Then define $g \in \Li(M)$ by
$$g = \sup_B (a_B - f_B)$$
where $a_B$ is the infimum of the essential range of $f|_B$, and the
supremum is taken in $L^\infty(M)$ (equivalently, as a limit
in the weak*-topology
of the unit ball of $\Li(M)$). Then $f \geq a_B - f_B$ for all $B$,
so $f \geq g|_A$. So if $f \neq g|_A$ then we must have
$f \geq g + \epsilon$ on some positive measure set $B \subset A$,
contradicting the definition of $a_B$ and the fact that $g \geq a_B - f_B$.
So $g$ has the desired properties.\hfill\hal
\bigskip

\noindent {\bf Theorem 29.} {\it Let $M$ be a measurable metric space and
$A \subset M$ a positive measure subset. Then $\X(A) = \chi_A\cdot \X(M)$.}
\medskip

\noindent {\it Proof.} The natural map from $\X(A)$ into $\chi_A\cdot \X(M)$
is isometric by Lemma 28. Conversely,
if $\phi \in \chi_A\cdot \X(M)$ and $f \in \Li(A)$, we can apply
$\phi$ to $f$ by first extending $f$ via Lemma 28; by Lemma 27 $\phi$
is insensitive to the extension. So $\X(A) = \chi_A\cdot\X(M)$.\hfill\hal
\bigskip

We now consider a special condition on $M$. We say a metric space is
``differentiable'' if, roughly speaking, its metric is captured by one-sided
derivations, as opposed to the two-sided derivations needed in general [69].
For geometric purposes differentiable spaces, or spaces which are locally
bi-Lipschitz equivalent to differentiable spaces, seem the most relevant.
The precise definition is as follows.
\bigskip

\noindent {\bf Definition 30.} The measurable metric space $M$ is
{\it differentiable}
if $\|df\| = L(f)$ for every $f \in \Li(M)$, where $d: \Li(M) \to \Omega(M)$
is the exterior derivative. Since $\|df\| \leq L(f)$ automatically, an
equivalent condition is
$$L(f) = \sup\{\|\phi f\|: \phi \in \X(M), \|\phi\| = 1\}$$
for all $f \in \Li(M)$.\hfill\hal
\bigskip

\noindent {\bf Theorem 31.} {\it Let $M = (M, \mu, \rho)$ be a complete
differentiable metric space and assume that every ball in $M$ with positive
radius has positive measure. Then $\rho$ is a path-metric.}
\medskip

\noindent {\it Proof.} Note that we are assuming $\rho$ is a pointwise
metric (justified, say, by Theorem 20). We use the term ``path-metric'' in
the sense of [9]: for each $x, y \in M$ the distance between $x$ and $y$ is
the infimum of $L(f)$ as $f$ ranges over all maps from $[0,1]$ into $M$
with $f(0) = x$ and $f(1) = y$.
\medskip

Let $\epsilon > 0$. Define a new metric $\rho_\epsilon$ by setting
$$\rho_\epsilon =
\inf \{\rho(x_0, x_1) + \rho(x_1, x_2) + \cdots + \rho(x_{n-1}, x_n)\},$$
where the infimum is taken over all finite sequences $x_0, \ldots, x_n$
such that $x_0 = x$, $x_n = y$, and $\rho(x_{i-1}, x_i) \leq
\epsilon$ ($1 \leq i \leq n$). It is easy to check that $\rho_\epsilon$
is a metric on $M$.
\medskip

Now fix $x, y \in M$ and define $f(z) = \rho_\epsilon(x,z)$. If
$\rho(z_1, z_2) \leq \epsilon$ then
$$|f(z_1) - f(z_2)| = |\rho_\epsilon(x,z_1) - \rho_\epsilon(x,z_2)|
\leq \rho_\epsilon(z_1, z_2) = \rho(z_1, z_2),$$
so $L(f|_A) \leq 1$ for any set $A \subset M$ with diameter at most $\epsilon$.
Then we can find $g \in \Li(M)$ such that $L(g) \leq 1$ (hence $\|dg\| \leq 1$)
and $g|_A = f|_A$ by Lemma 28. Thus $(df)|_A = (dg)|_A$ by Lemma 27;
as $M$ can be covered by positive measure balls of diameter at most $\epsilon$,
we conclude that $\|df\| \leq 1$. This implies
$L(f) \leq 1$ by differentiability, so
$$\rho_\epsilon(x,y) = |f(x) - f(y)| \leq \rho(x,y).$$
As the reverse inequality is automatic, we have $\rho(x,y) =
\rho_\epsilon(x,y)$, and equality for all $\epsilon$ plus completeness
implies that $\rho$ is a path-metric ([30], Th\'eor\`eme 1.8).\hfill\hal
\bigskip

The converse of Theorem 31 is false; there exist path-metric spaces $M$
for which $\X(M) = 0$ (see Theorem 40).
\medskip

An easy nontrivial example of a
differentiable space in which every ball of finite radius has measure zero
is $\R^2$ with Lebesgue measure and metric
$$\rho((x,y), (x', y')) = \cases{|x - x'|& if $y = y'$\cr
\infty& if $y \neq y'$\cr}.$$
Clearly, this space will remain differentiable if its metric
is modified on a single horizontal line, but this can be done in such a way
that $\rho$ is no longer a path-metric.
\bigskip
\bigskip

\noindent {\underbar{\bf 5. Examples.}}
\bigskip

We now determine $\X(M)$ for various spaces $M$.
Theorem 29 will be used repeatedly; it implies that our analysis need only
be done locally. After restriction to a manageable subset of $M$, the main
tools are then Theorems 26 and 16.
\bigskip

\noindent {\bf A. Atomic measures and Stone spaces}
\bigskip

We begin with some examples of metric spaces which admit no metric derivations,
and hence are zero-dimensional in the sense that $\X(M) = 0$. The first
result shows that the measure really is an essential ingredient in the
construction of $\X(M)$ (cf.\ Corollary 7 of [69]).
\bigskip

\noindent {\bf Proposition 32.} {\it Let $M$ be a metric space equipped with
an atomic measure. Then $\X(M) = 0$.}
\medskip

\noindent {\it Proof.} Let $\delta \in \X(M)$. Also
let $x \in M$ and let $I = \{f \in \Li(M): f(x) = 0\}$.
Then $I$ is a weak*-closed ideal of $\Li(M)$, so we have $\delta(I)
\subset I$ just as in the proof of Lemma 27. Now letting $1$ denote
the function which is constantly 1, we have $\delta(1) = 0$
since $\delta(1\cdot 1) = 1\cdot\delta(1) + \delta(1)\cdot 1 = 2\delta(1)$.
Thus $f - f(x)\cdot 1 \in I$ implies
$$\delta(f) = \delta(f - f(x)\cdot 1) \in I,$$
i.e.\ $(\delta f)(x) = 0$. This holds for all $x \in M$, so
$\delta f = 0$. This shows that $\X(M) = 0$.\hfill\hal
\bigskip

As a technical note, we should point out that our restriction to
$\sigma$-finite measure spaces severely restricts the scope of Proposition 32.
However, any $l^\infty(X)$ is obviously
finitely decomposable, regardless of the cardinality of $X$,
and Proposition 32 in fact remains true at that level of generality. (Indeed,
there is no real difficulty in doing everything up to this point with
finitely decomposable measures. But it seems preferrable to work with the
more familiar property of $\sigma$-finiteness, since the only apparent
drawback in doing so
is the exclusion of uninteresting examples like $l^\infty(X)$.)
\medskip

The next result shows that metric spaces with a certain disconnectedness
property have no nonzero metric derivations. The condition implies that
$M$ is totally disconnected, but not every totally disconnected space
satisfies it; in fact, any totally disconnected subset of $\R^n$ with
positive measure will have nontrivial $\X(M)$ (Theorem 37).
\bigskip

\noindent {\bf Proposition 33.} {\it Let $M$ be a measurable metric space
and suppose the simple Lipschitz functions
are weak*-dense in $\Li(M)$. Then $\X(M) = 0$.}
\medskip

\noindent {\it Proof.} Let $\delta \in \X(M)$ and let
$f = \sum a_n \chi_{A_n} \in \Li(M)$ be a simple function with the sets
$A_n$ disjoint. Then $f =
a_n\cdot 1$ on $A_n$, so $\delta(f)|_{A_n} = 0$ by Lemma 27; as this
is true for all $n$, we have $\delta(f) = 0$. Density of the simple
functions then implies that $\delta = 0$.\hfill\hal
\bigskip

\noindent {\bf Corollary 34.} {\it Let $M$ be a measurable metric space
which is uniformly discrete in the sense that there exists $\epsilon > 0$
such that every $A, B \subset M$ satisfy $\rho(A, B) = 0$ or $\rho(A, B)
\geq \epsilon$. Then $\X(M) = 0$.}
\bigskip

\noindent {\it Proof.} The uniform discreteness condition implies that
$\Li(M) \cong L^\infty(M)$. So the simple functions are actually norm-dense
in $\Li(M)$.\hfill\hal
\bigskip

\noindent {\bf Corollary 35.} {\it Let $K$ be the middle-thirds Cantor set,
equipped with any $\sigma$-finite Borel measure and with metric inherited
from $\R$. Then $\X(K) = 0$.}
\medskip

\noindent {\it Proof.} The simple functions are weak*-dense in $\Li(K)$
by Theorem 16 with $k = 3$.\hfill\hal
\bigskip

\noindent {\bf B. Lipschitz manifolds}
\bigskip

Next we consider the case of Lipschitz manifolds ([17], [43], [52]).
Note that this class includes all $C^1$-Riemannian manifolds.
\bigskip

\noindent {\bf Theorem 36.} {\it Let $M$ be a Lipschitz manifold, equipped
with Lebesgue measure class. Then $\X(M)$, $\Omega(M)$, and $d$ have
their usual meanings (as in [27]).}
\medskip

\noindent {\it Proof.} We check that every bounded measurable
vector field (in the usual sense)
gives rise to a metric derivation, and vice versa. By Theorem 29
it suffices to consider the case that $M$ is a region in $\R^n$.
\medskip

First consider the vector field $\partial/\partial x_k$. Every Lipschitz
function is almost everywhere differentiable in any coordinate direction,
with derivative in $L^\infty(M)$ ([22], Theorem 3.1.6);
and if $f_i \to f$ boundedly weak* in $\Li(M)$ then
$(\partial f_i/\partial x_k)$ is bounded and
$$\int_M \big({{\partial f_i}\over{\partial x_k}}\big) g
= -\int_M f_i{{\partial g}\over{\partial x_k}}
\to -\int_M f{{\partial g}\over{\partial x_k}}
= \int_M \big({{\partial f}\over{\partial x_k}}\big) g$$
for every $g \in C^\infty(M)$ supported on the interior of $M$,
which shows that $\partial f_i/\partial x_k \to \partial f/\partial x_k$
weak* in $L^\infty(M)$. So $\delta(f) = \partial f/\partial x_k$ is
a metric derivation.
\medskip

Since the metric derivations constitute an $L^\infty(M)$-module, it follows
that the vector field
$\sum_1^n f_k(\partial/\partial x_k)$ gives rise to a metric derivation for
any $f_1, \ldots, f_n \in L^\infty(M)$. So every bounded measurable vector
field gives rise to a metric derivation.
\medskip

Conversely, any metric derivation $\delta$ is determined by its values
$f_k = \delta(x_k)$ on the coordinate functions since these functions
generate $\Li(M)$ by Theorem 16. So there are no metric derivations
besides those which arise from bounded measurable vector fields.\hfill\hal
\bigskip
\bigskip

\noindent {\bf C. Rectifiable sets}
\bigskip

We consider $(H^m, m)$ rectifiable and $H^m$ measurable subsets of $\R^n$
in the sense of ([22], 3.2.14). Here $H^m$ is $m$-dimensional Hausdorff measure.
\bigskip

\noindent {\bf Theorem 37.} {\it Let $M$ be an $(H^m, m)$ rectifiable and
$H^m$ measurable subset of $\R^n$. Then $\X(M)$ is naturally identified
with the module of bounded measurable sections of approximate tangent
spaces ([22], 3.2.16 and 3.2.19).}
\medskip

\noindent {\it Proof.} By ([22], Lemma 3.2.18) $M$ can be decomposed into
a countable union of sets each of which is bi-Lipschitz equivalent to
a positive measure subset of $\R^m$. So by Theorem 29 we can reduce to the
case that $M$ is a positive measure subset of $\R^m$. In this case the
approximate tangent space is $\R^m$ at almost every point ([22], Theorem
3.2.19), and the identification of metric derivations with bounded measurable
vector fields then follows from Theorem 36 and Theorem 29.\hfill\hal
\bigskip

\noindent {\bf D. Sub-Riemannian metrics}
\bigskip

Let $M$ be an $n$-dimensional Riemannian manifold and let $B$ be a smooth
$k$-dimensional subbundle of the tangent bundle $TM$. Define the length of
any smooth path $\gamma: [0,1] \to M$ to be
$l(\gamma) = \int_0^1 \|\gamma'\| dt$, and define a metric $\rho$ on
$M$ by setting
$$\rho(x,y) = \inf\{l(\gamma): \gamma(0) = x,\,
\gamma(1) = y,\hbox{ and }\gamma'(t) \in B\hbox{ for all }t \in [0,1]\}.$$
Let $\rho'$ denote the usual Riemannian metric, defined by the same
infimum taken over all smooth paths from $x$ to $y$.
Equip $M$ with the usual Lebesgue measure class.
\bigskip

\noindent {\bf Theorem 38.} {\it $\X(M)$ is naturally identified with the
bounded measurable sections of $B$.}
\medskip

\noindent {\it Proof.} By Theorem 29 it is sufficient to consider a small
neighborhood of any point $x \in M$. Fix $x$ and let $\alpha$ be a
diffeomorphism between a neighborhood of $x$ and an open set in $\R^n$
such that $\alpha(x) = 0$.
Let $v_1, \ldots, v_n \in T_xM$ be orthonormal vectors, each of which
either belongs to or is orthogonal to $B_x$, and let $X_1, \ldots, X_n$ be
smooth vector fields on $M$ such that $X_i(x) = v_i$. By projecting
onto $B$ and $B^\perp$ and then
orthonormalizing, we may assume that in a neighborhood of $x$ the $X_i$ are
orthonormal ($1 \leq i \leq n$) and $X_1, \ldots, X_k$ span $B$.
Then, by composing $\alpha$ with an invertible linear transformation on
$\R^n$, we may assume that the vectors $\alpha_* v_i$ are orthonormal with
respect to the Euclidean metric on $\R^n$ at 0. On a small enough
neighborhood of $0$ the vector fields $\alpha_* X_i$ are then
nearly orthonormal with respect to the Euclidean metric, and using a
partition of unity we may (1) find a metric on $\R^n$ which agrees with
the Euclidean metric outside a neighborhood of
$0$ and makes $\alpha$ isometric near $0$ and (2) extend the
$\alpha_* X_i$ to linearly independent vector fields on all of $\R^n$.
Finally, applying Gramm-Schmidt we may take the $X_i$ to
be orthonormal with respect to the metric just introduced.
Thus, in what follows we will
assume that $M = \R^n$ with a metric which is Euclidean outside a
bounded region, and that there are globally defined orthonormal vector fields
$X_1, \ldots, X_n$ the first $k$ of which span $B$. The reduction given
in this paragraph is due to Renato Feres [24].
\medskip

Now for $1 \leq i \leq n$ let $T^i_t$ ($t \geq 0$) be the flow generated by
$X_i$, and define $\alpha^i_t: L^\infty(M) \to L^\infty(M)$ by
$\alpha^i_t(f) = f\circ T^i_t$. Then for each $i$, $(\alpha^i_t)$ is a
strongly continuous one-parameter group of automorphisms of $L^\infty(M)$.
Let $\delta_i: L^\infty(M) \to L^\infty(M)$ be its infinitesimal generator,
defined by $\delta(f) = \lim_{t\to 0} (f - \alpha^i_t(f))/t$, with domain
consisting
of all $f \in L^\infty(M)$ for which the limit exists in the weak* sense.
\medskip

Suppose $1 \leq i \leq k$. If $x \in M$ then $\gamma(t) = T^i_t(x)$ is a path
whose tangent vector lies in $B$ and has norm equal to one everywhere. Thus
$\rho(x, T^i_t(x)) \leq |t|$ for all $t$. So for any $f \in \Li(M)$ we have
$$|f(x) - \alpha^i_t(f)(x)| = |f(x) - f(T^i_t(x))| \leq |t|L(f);$$
this shows that $\|f - \alpha^i_t(f)\| \leq |t|L(f)$,
and ([10], Proposition 3.1.23) then implies that $f$ belongs to the domain
of $\delta_i$. In fact, the above inequality shows that $\|\delta_i(f)\|
\leq L(f)$, so that $\delta_i$ is nonexpansive when regarded as a map from
$\Li(M)$ to $L^\infty(M)$. Furthermore, $\delta_i$ is a W*-derivation by
([10], Proposition 3.1.6) so its restriction to $\Li(M)$ is a metric
derivation. Now for any $f_1, \ldots, f_k \in L^\infty(M)$, we have
$\sum_1^k f_i\delta_i \in \X(M)$, so every bounded measurable section
of $B$ defines a metric derivation.
\medskip

To prove the converse, we invoke Theorem 26. Let $E$ be the set of metric
derivations of the form $\sum_1^k f_i\delta_i$ described above. It is
isometrically isomorphic as an $L^\infty(M)$-normed module to
$L^\infty(M, \R^n)$ (giving $\R^n$ the
Euclidean norm) and hence is reflexive.
For the other hypothesis, let $f \in \Li(M)$. Convolving $f$ with a
$C^\infty$ approximate unit of $L^1(\R^n)$ produces a sequence of smooth
functions, bounded in Lipschitz norm, which converge to $f$ weak* in
$L^\infty(M)$. This shows that the smooth functions are weak*-dense in
$\Li(M)$, so we may assume $f$ is smooth.
\medskip

In the notation of Theorem 26, we must show that $|d_E f| \geq |df|$.
(The reverse inequality is automatic.) Since $\|df|_A\| \leq L(f|_A)$ for
any positive measure set $A$, it will suffice to find, for any
$x \in M$ and $\epsilon > 0$, a neighborhood $A$ of $x$ such that
$|d_E f| \geq L(f|_A) - \epsilon$ on $A$.
\medskip

Let $a = (\sum_1^n (\delta_i f(x))^2)^{1/2}$ and $b =
(\sum_1^k (\delta_i f(x))^2)^{1/2}$. Since $f$ is smooth, we may
find $r > 0$ such that $b - \epsilon/2 \leq |d_E f| \leq b + \epsilon/2$
and $|df| \leq a + \epsilon$ on the $\rho'$-ball of radius $r$ about $x$.
(The bound on $|df|$ is due to the fact that every metric derivation of
$\Li(M, \rho)$ is also a metric derivation of $\Li(M, \rho')$, hence is
a linear combination of partial derivatives by Theorem 36.)
Let $s = r(b+\epsilon/2)/3(a+\epsilon)$, and let
$A$ be the $\rho'$-ball of radius $s$ about $x$.
\medskip

Let $y, z \in A$; we must show that $|f(y) - f(z)|/\rho(y,z) \leq
b + \epsilon/2$. Now $\rho'(y,z) \leq 2s$, so
$$|f(y) - f(z)| \leq (a+\epsilon)\rho'(y,z) \leq 2r(b+\epsilon/2)/3.$$
This shows that the desired inequality holds if $\rho(y,z) \geq 2r/3$.
Otherwise, let $\gamma: [0,1] \to M$ be a constant velocity path from
$y$ to $z$ which is everywhere tangent to $B$ and whose total length is
exactly $\rho(y,z) < 2r/3$; this exists by ([9], Lemma 2.1.2). Since
$\rho'(y,x) \leq s < r/3$ it follows that $\gamma$ lies entirely within
the $\rho'$-ball of radius $r$ about $x$, so that $|d_E f|(\gamma(t))
\leq b + \epsilon/2$ for any $t \in [0,1]$. Thus
$$|f(y) - f(z)| \leq l(\gamma)\sup_{t \in [0,1]} |d_E f|(\gamma(t))
\leq \rho(y,z)(b + \epsilon/2),$$
as desired.\hfill\hal
\bigskip

Theorem 38 is no longer true if we allow the dimension of $B$ to vary;
this is illustrated by the space $M_1$ treated in Theorem 55 (see the
comment following that theorem).
\bigskip

\noindent {\bf E. The Sierpinski carpet}
\bigskip

Let $S$ be the Sierpinski carpet obtained from the unit square by iterating
the process of removing the middle ninth sub-square. That is, $S$ is
the set of points $x = (x_1, x_2) \in [0,1]^2$ such that 
for no $n \in \N$ and $1 \leq k, l \leq 3^{n-1}$ is it the case that
$${{3k-2}\over{3^n}} < x_1 < {{3k-1}\over{3^n}}\qquad{\rm and}\qquad
{{3l-2}\over{3^n}} < x_2 < {{3l-1}\over{3^n}}.$$
We give $S$ normalized Hausdorff measure $\mu$; this means that if
$$S_{k,l,n} = S \cap ({{k-1}\over{3^n}}, {{k}\over{3^n}})\times
({{l-1}\over{3^n}}, {{l}\over{3^n}})$$
is nonempty then $\mu(S_{k,l,n}) = 8^{-n}$.
\bigskip

\noindent {\bf Lemma 39.} {\it Let $f \in L^\infty(S)$ and let
$a$ belong to the essential range of $f$. For any $\epsilon_1, \epsilon_2 > 0$
there exist $k, l, n$ so that
$$\mu(S_{k,l,n} \cap f^{-1}([a - \epsilon_1, a+ \epsilon_1])) \geq
8^{-n}(1-\epsilon_2).$$}

\noindent {\it Proof.} The usual proof of the Lebesgue differentiation
theorem can be adapted to show that for almost every $x \in S$ we have
$$\lim 8^n\int_{S_{k,l,n}} |f(y) - f(x)|dy = 0,$$
for any sequence of squares $S_{k,l,n}$ each of which contains $x$ and
for which $n \to \infty$. For instance, the argument in ([25], $\S$ 3.4)
can be carried over verbatim, replacing $\R^n$ with $S$ and ``open ball''
with ``the closure of some $S_{k,l,n}$.''
\medskip

Since this holds for almost every $x \in S$, it must hold for some $x_0 \in
f^{-1}((a - \epsilon_1, a+\epsilon_1))$ by the definition of essential
range. Choose $\epsilon$ such that
$$[f(x_0) - \epsilon, f(x_0) + \epsilon]
\subset [a - \epsilon_1, a+\epsilon_1];$$
then find a square $S_{k,l,n}$ which contains $x_0$ such that
$$8^n\int_{S_{k,l,n}} |f(y) - f(x_0)| dy \leq \epsilon\epsilon_2.$$
It follows that
$$8^n\mu(S_{k,l,n} \cap f^{-1}([f(x_0) - \epsilon, f(x_0) + \epsilon])
\geq 1-\epsilon_2,$$
which is enough.\hfill\hal
\bigskip

\noindent {\bf Theorem 40.} {\it $\X(S) = 0$.}
\medskip

\noindent {\it Proof.} Let $\delta \in \X(S)$. It will be enough to show
that $\delta(f) = 0$ when $f$ is either of the two coordinate functions,
$f(x,y) = x$ or $f(x,y) = y$, since these generate Lip(S) by Theorem 16.
The arguments in the two cases are the same, so take
the first case and suppose $\delta(f) \neq 0$. Since this
implies $(a\delta)(f) \neq 0$ for any nonzero constant $a$, we can suppose
without loss of generality that $\|\delta(f)\| = {\rm ess}\sup \delta(f) = 1$.
\medskip

Define a sequence of piecewise-linear functions $f_m \in C(S)$ by letting
$f_m(0,y) = 0$ and requiring
$${{\partial f_m}\over{\partial x}}(x,y) =
\cases{1&if $3k/3^m < x < (3k+1)/3^m$\cr
-1&if $(3k+1)/3^m < x < (3k+2)/3^m$\cr
0&if $(3k+2)/3^m < x < (3k+3)/3^m$\cr}.$$
Then $f_m \to 0$ weak* in $\Li(S)$. Also
$f_m - f$ is constant on the left $3/8$ of each $S_{k,l,m}$; $f_m + f$ is
constant on the middle $2/8$ of each $S_{k,l,m}$; and
$f_m$ is zero on the right $3/8$ of each $S_{k,l,m}$.
\medskip

Choose $\epsilon_1, \epsilon_2 > 0$ such that $c = 1/8 - 3\epsilon_1/8 -
2\epsilon_2 > 0$, and apply the lemma to find a square $S_{k,l,n}$ such that
$$\mu^{-1}(S_{k,l,n} \cap (\delta f)^{-1}([1 - \epsilon_1, 1 + \epsilon_1]))
\geq 8^{-n}(1-\epsilon_2).$$
Now let $m \geq n$ and
consider $8^n\int_{S_{k,l,n}} \delta(f_m)$. This integral is zero on the
$3/8$ of $S_{k,l,n}$ where $f_m$ is zero. On the
$2/8$ where $f_m + f$ is constant, we have $|\delta(f_m)| = |\delta(f)|
\leq 1$ and so the integral is not below $-2/8$. Of the remainder,
$\delta(f_m) \geq 1-\epsilon_1$ on a set of measure at least
$8^{-n}(3/8 - \epsilon_2)$ and $|\delta(f_m)| = |\delta(f)| \leq 1$
elsewhere, so the integral here is at least
$(3/8 - \epsilon_2)(1-\epsilon_1) - \epsilon_2$. All together, we have
$$\int_{S_{k,l,n}} \delta(f_m) \geq
8^{-n}\big((3/8 - \epsilon_2)(1 - \epsilon_1) - \epsilon_2 - 2/8\big)
> c/8^n.$$
This holds for every $m \geq n$,
so $\int \delta(f_m)\chi_{S_{k,l,n}}$ does not go to
zero as $m \to \infty$, contradicting weak*-continuity of $\delta$. This shows
that the assumption $\delta(f) \neq 0$ is impossible. So $\delta$ vanishes
on the coordinate functions, hence $\delta = 0$.\hfill\hal
\bigskip

A similar argument shows that the Sierpinski gasket (obtained by iterating
the process of removing the middle fourth sub-triangle from an equilateral
triangle) supports no nonzero metric derivations.
I have not tried to systematically extend the reasoning in Theorem 40 to
other fractals, but it seems likely that the same sort of argument would
apply to many fractal shapes with non-integral Hausdorff dimension. 
\medskip

The Sierpinski carpet is the closure of a sequence of finite graphs
$G_n$. Namely, $G_1$ is the boundary of the unit square $[0,1]^2$ and
$G_{n+1} = \bigcup_{k=1}^8 (G_n + v_k)/3$ where $v_1, \ldots, v_8$
are the vectors $(0,0)$, $(1, 0)$, $(2, 0)$, $(0, 1)$, $(2, 1)$,
$(0, 2)$, $(1, 2)$, $(2, 2)$. It may be worth noting that
a reasonable one-dimensional differentiable structure on $S$ can
be obtained by setting aside Hausdorff measure and instead assigning
zero measure to $S - \bigcup G_n$ and using one-dimensional Lebesgue
measure on each $G_n$.
\bigskip

\noindent {\bf F. Hilbert cubes}
\bigskip

Fix $1 < p < \infty$ and a sequence $(a_n) \in l^p(\N)$ with $a_n > 0$
for all $n$. Let $M_p$ be the cartesian product
$M_p = \prod [0, a_n]$ with metric inherited from
$l^p(\N)$. Also, give $M_p$ the product of normalized Lebesgue measure on
each factor. Let $q$ be the conjugate exponent to $p$.
\medskip

We say that a sequence $(f_n) \subset L^\infty(M_p)$ is {\it weakly
$p$-summable} if the partial sums $\sum_1^N |f_n|^p$ are uniformly
bounded in $L^\infty(M_p)$. We then define $|(f_n)|^p = \sum |f_n|^p$ to be the
supremum in $L^\infty(M_p)$ of the partial sums, and denote the
$L^\infty(M_p)$-normed module of all weakly $p$-summable sequences by
$l^p(L^\infty(M_p))$.
\bigskip

\noindent {\bf Lemma 41.} {\it For any $f \in \Li(M_p)$, the sequence
$(\partial f/\partial x_n)$ is weakly $q$-summable.}
\medskip

\noindent {\it Proof.} Note that $\partial f/\partial x_n$ exists almost
everywhere on $M$, by Rademacher's theorem ([22], Theorem 3.1.6) plus
Fubini's theorem. We will show that
$\sum_1^N |\partial f/\partial x_n|^q \leq L(f)$ almost everywhere, for
any $N \in \N$. By
Fubini's theorem it suffices to show this for Lipschitz functions on
$M_p^N = \prod_1^N [0,a_n]$.
\medskip

Let $f \in \Li(M_p^N)$ and suppose $f$ is differentiable at the point
$x = (x_1, \ldots, x_N) \in M_p^N$. Then for any $b = (b_1, \ldots, b_N)
\in \R^N$ and $\epsilon > 0$ there exists $r > 0$ such that
$$|f(x + r b) - f(x)| \geq
(1-\epsilon)r|\sum b_n(\partial f/\partial x_n)(x)|.$$
But also
$$|f(x + r b) - f(x)| \leq r L(f)\sum |b_n|^p,$$
so taking $\epsilon \to 0$ we have
$$|\sum b_n(\partial f/\partial x_n)(x)| \leq L(f) \sum |b_n|^p.$$
As this holds for any $b \in \R^N$, we conclude that
$\sum_1^N |\partial f/\partial x_n(x)|^q \leq L(f)$ at every point $x$ where
$f$ is differentiable. Since any Lipschitz function is differentiable almost
everywhere, we are done.\hfill\hal
\bigskip

\noindent {\bf Lemma 42.} {\it Let $(f_n) \in l^p(L^\infty(M_p))$
and let $\epsilon_1, \epsilon_2 > 0$. Then there
exists $N \in \N$ and a subset $A_N \subset M_p$ such that $\mu(M_p - A_N)
\leq \epsilon_1$ and $\sum_{N+1}^\infty |f_n|^p \leq \epsilon_2$ almost
everywhere on $A_N$.}
\medskip

\noindent {\it Proof.} Fix a Borel representative of each $f_n$. For
each $N$ let
$$A_N = \{x \in M_p: \sum_{N+1}^\infty |f_n(x)|^p \leq \epsilon_2\}.$$
Then the $A_N$ are nested and $\bigcup A_N = M_p$; since $M_p$ has finite measure
this implies that $\mu(M_p - A_N) \leq \epsilon_1$ for some $N$.\hfill\hal
\bigskip

\noindent {\bf Lemma 43.} {\it $l^p(L^\infty(M_p))' \cong l^q(L^\infty(M_p))$.}
\medskip

\noindent {\it Proof.} The product of a sequence in $l^p(L^\infty(M_p))$
and a sequence in $l^q(L^\infty(M_p))$ is weakly $1$-summable, hence
converges almost everywhere, hence is bounded and converges weak* in
$L^\infty(M_p)$. This shows that $l^q(L^\infty(M_p))$ naturally embeds in
$l^p(L^\infty(M_p))'$, and the embedding is clearly isometric.
\medskip

Let $\Phi \in l^p(L^\infty(M_p))'$ and for each $m \in \N$ let $g_m =
\Phi(\delta_{m,n}\cdot 1_{M_p})$ where $\delta_{m,n}\cdot 1_{M_p}$ is
the sequence whose $m$th term is the constant function $1_{M_p}
\in L^\infty(M_p)$ and whose other terms are all zero.
We claim that $(g_m) \in l^q(L^\infty(M_p))$.
If not, then for every $C > 0$ there exists $N \in \N$ and a positive
measure set $A \subset M_p$ such that $\sum_1^N |g_m|^q \geq C$ on $A$.
By shrinking $A$ we may assume that each $g_m$ ($1 \leq m \leq N$)
varies by less than $N^{-1/q}$ on $A$. For $1 \leq m \leq N$ let
$a_m$ be a constant such that $|g_m - a_m| \leq N^{-1/q}$ on $A$, so that
$\sum_1^N |a_m|^q \geq C - 1$. Then find a sequence $(b_m)$ with $l^p$-norm
one such that $\sum_1^N a_mb_m \geq C - 1$; let $f_m = b_m\cdot 1_{M_p}$;
and observe that $f = (f_m) \in l^p(L^\infty(M_p))$ and $|(f_m)|_p = 1_{M_p}$
but
$$\Phi(f) = \sum_1^N f_mg_m \geq C - 2$$
on $A$. This contradicts boundedness of $\Phi$ and establishes that
$(g_m) \in l^q(L^\infty(M_p))$.
\medskip

Finally, we must show that $\Phi$ is given by summation against $(g_m)$.
This is clearly true on any $p$-summable sequence which has only finitely
many nonzero terms. Now suppose $\Phi$ annihilates every such
sequence. Then for any $f = (f_n) \in l^p(L^\infty(M_p))$ and any
$\epsilon_1, \epsilon_2 > 0$ we can apply Lemma 42
to get a set $A_N$ of measure at least $1 - \epsilon_1$
such that $(\chi_{A_N}f_n)$ is within $\epsilon_2$
in norm from a sequence with only finitely many nonzero terms. Thus
$|\Phi(f)| \leq \epsilon_2\|\Phi\|$ on $A_N$. Taking $\epsilon_1,
\epsilon_2 \to 0$ shows that $\Phi(f) = 0$. This shows that every element
of $l^p(L^\infty(M_p))'$ is determined by its values on finite sequences, which completes
the proof.\hfill\hal
\bigskip

\noindent {\bf Theorem 44.} {\it $\X(M_p) \cong l^p(L^\infty(M_p))$.}
\medskip

\noindent {\it Proof.} Let $g = (g_n) \in l^p(L^\infty(M_p))$.
For any $f \in \Li(M_p)$ the series
$$\delta_g(f) = \sum g_n(\partial f/\partial x_n)$$
converges weak* in $L^\infty(M_p)$ by Lemmas 41 and 43.
\medskip

We claim that $\delta_g$ is a metric derivation. Linearity and the derivation
identity are easy, as is the inequality $\|\delta_g\| \leq \|g\|$.
To check weak*-continuity, suppose $f_i \to f$ weak*
in the unit ball of $\Li(M_p)$. By taking a subnet, we may assume that
$\delta_g(f_i)$ converges weak* to some $h \in L^\infty(M_p)$; we must show
that $h = \delta_g(f)$. Given $\epsilon_1, \epsilon_2 > 0$ find a set
$A_N$ as in Lemma 42 for the sequence $(g_n)$.
Let $g^N_n = g_n$ if $n \leq N$ and 0 otherwise.
Then $\delta_{g^N}$ is weak*-continuous since it is a finite linear combination
of partial derivatives, each of which is weak*-continuous by the Theorem 36
plus Fubini's theorem. Also
$$\delta_g(f_i) - \delta_g(f)
= (\delta_g(f_i) - \delta_{g^N}(f_i)) + \delta_{g^N}(f_i - f)
+ (\delta_{g^N}(f) - \delta_g(f)),$$
and on $A_N$ the first and third terms are each at most $\epsilon_2$ in
absolute value. Taking the limit,
weak*-continuity of $\delta_{g^N}$ implies that $|h - \delta_g(f)| \leq
2\epsilon_2$ on $A_N$. Then taking $\epsilon_1, \epsilon_2 \to 0$ establishes
that $h = \delta_g(f)$, and we conclude that $\delta_g$ is a metric
derivation. Thus every $p$-summable sequence
in $L^\infty(M_p)$ gives rise to an element of $\X(M_p)$.
\medskip

We invoke Theorem 26 to prove the converse. Let $E$ be the set of all
metric derivations which arise from $p$-summable sequences. It is reflexive
by Lemma 43. For the other hypothesis, consider the
Lipschitz functions on $M_p$ of the form $f \circ \pi^N$ for some
$N \in \N$ and $f \in \Li(M_p^N)$, where $M_p^N = \prod_1^N [0, a_n]$
and $\pi^N: M_p \to M_p^N$ is the natural projection. By Theorem 16 these
functions are weak*-dense in $\Li(M_p)$, so it will suffice to verify that
$|d(f\circ \pi^N)| = |d_E(f\circ \pi^N)|$ for $f \in \Li(M_p^N)$.
\medskip

Let $d^N$ be the exterior derivative on $M_p^N$. By Theorem 36 we know
that the coordinate partial derivatives generate $\X(M_p^N)$, so
$$|d_E(f\circ \pi^N)|(z) = |d^N f|(\pi^N(z))$$
for almost all $z \in M_p$. Now let $\delta \in \X(M_p)$,
$\|\delta\| \leq 1$. For any norm-one
$L^\infty(M_p^N)$-module map $\Phi: L^\infty(M_p) \to L^\infty(M_p^N)$,
the map $f \mapsto \Phi(\delta(f\circ \pi^N))$ is a metric derivation
on $M_p^N$, hence
$$|\Phi(\delta(f\circ\pi^N))| \leq |d^N f|$$
on $M_p^N$. As this is true for all $\Phi$, Corollary 6 implies that
$$|\delta(f\circ \pi^N)|(z) \leq |d^N f|(\pi^N(z))$$
almost everywhere; so we conclude that
$$|d(f\circ \pi^N)|(z) \leq |d^N f|(\pi^N(z)) = |d_E(f\circ \pi^N)|(z)$$
for almost all $z \in M_p$. As the reverse inequality is automatic, this
verifies the second hypothesis of Theorem 26 and therefore completes the
proof.\hfill\hal
\bigskip

The cases $p = 1$ and $p = \infty$ are less transparent, but they are still
of interest because they allow us to falsify some natural conjectures. Thus,
fix $(a_n) \in l^1(\N)$ with $a_n > 0$ for all $n$ and
let $M_1 = \prod [0,a_n]$. Also let $(b_n) \in c_0(\N)$ with $b_n > 0$
for all $n$ and let $M_0 = \prod [0, b_n]$. Give $M_1$ the $l^1$
metric and $M_0$ the $c_0$ metric, and endow both with the product of
normalized Lebesgue measure on each factor.
\bigskip

\noindent {\bf Proposition 45.} {\it (a). $\X(M_0)$ is not weak*-closed
in $B(\Li(M_0), L^\infty(M_0))$ (the space of bounded linear
maps from $\Li(M_0)$ into $L^\infty(M_0)$).
\medskip

(b). There is a metric derivation $\delta \in \X(M_1)$ which is not a
W*-derivation.}
\medskip

\noindent {\it Proof.} (a). For each $n \in \N$ define $\delta_n \in
\X(M_0)$ by
$$\delta_n(f) = \sum_{k = 1}^n (\partial f/\partial x_k),$$
and let $\delta$ be a weak*-cluster point of the (bounded)
sequence $(\delta_n)$ in the dual space
$B(\Li(M_0), L^\infty(M_0))$.
\medskip

Let $f_m \in \Li(M_0)$ be the $m$th coordinate function, $f_m(x) = x_m$.
Then $f_m \to 0$ uniformly, hence weak* in $\Li(M_0)$. However
$\delta(f_m) = \lim \delta_n(f_m) = 1_{M_0}$ for all $m$,
so $\delta(f_m)$ does not
converge to zero weak*, and hence $\delta$ cannot be a metric derivation.
\medskip

(b). Define $\delta: \Li(M_1) \to L^\infty(M_1)$ by
$$\delta(f) = \sum_{k=1}^\infty a_k(\partial f/\partial x_k).$$
Then $\delta \in \X(M_1)$ because $\X(M_1)$ is a Banach space and the
sequence $(a_n)$ is summable. For each $n \in \N$ define $f_n(x) = x_n/a_n$.
Then $f_n \to (1/2)\cdot 1_{M_1}$ weak* in $L^\infty(M_1)$, while
$\delta(f_n) = 1_{M_1}$ for all $n$. Since $\delta((1/2)\cdot 1_{M_1}) = 0$,
this shows that $\delta$ is not a W*-derivation.\hfill\hal
\bigskip

\noindent {\bf G. Banach manifolds}
\bigskip

For our purposes a Banach manifold is a metric space $M$ which is locally
bi-Lipschitz equivalent to the unit ball of some fixed Banach space.
We appeal to Theorem 29 to reduce to the case where $M = F$ is itself a
Banach space. This ignores the issue that there is in general no canonical
choice of measure or measure class on (the unit ball of)
an infinite-dimensional Banach space, so measure-theoretic complications
may arise when one tries to build a manifold by patching
local neighborhoods together.
\medskip

Let $F$ be a separable
reflexive Banach space equipped with a $\sigma$-finite Borel
measure $\mu$. We assume that there is a dense subspace $F_0 \subset F$
with the property that $\mu$ and its translation $\mu_v$ by $v$ are mutually
absolutely continuous for any $v \in F_0$.
\bigskip

\noindent {\bf Lemma 46.} {\it Let $X$ be a $\sigma$-finite measure space.
Then $L^\infty(X, F)$ is an $L^\infty(X)$-normed module and
$L^\infty(X, F)' \cong L^\infty(X, F')$.}
\medskip

\noindent {\it Proof.} By $L^\infty(X, F)$ we mean the space of bounded
measurable functions from $X$ into $F$, modulo functions which
vanish off of a null set. Verification of the
$L^\infty(X)$-normed module property given in Theorem 2 (c) is easy. For the
second assertion,
if $\Phi \in L^\infty(X, F')$ and $\phi \in L^\infty(X, F)$ then
the map $x \mapsto P(\Phi(x), \phi(x))$ is a bounded measurable function
on $X$, where $P: F'\times F \to \R$ is the natural pairing. From here it is
straightforward to check that $L^\infty(X, F')$ embeds isometrically in
$L^\infty(X, F)'$. For instance, this can be done by showing that the
simple functions in $L^\infty(X, F')$ are norm-dense, and checking isometry
on them.
\medskip

Conversely, let $\Phi \in L^\infty(X, F)'$. Let $S \subset F$ be a countable
dense subset. For each $v \in S$ fix a Borel version $f_v$ of
$\Phi(v\cdot 1_X)$. Then for any finite subset $S_0$ of $S$ the map
$$\sum_{v \in S_0} a_v v \mapsto \sum_{v \in S_0} a_vf_v(x)$$
is a bounded linear functional, of norm at most $\|\Phi\|$, on the linear span
of $S_0$, for almost every $x \in X$.
Thus there is a set $X_0 \subset X$ of full measure such
that for all $x \in X_0$ the map $\sum a_iv_i \mapsto \sum a_if_{v_i}(x)$
is a bounded linear functional on the unclosed span of $S$. So
for each $x \in X_0$ there exists a unique element $\Psi(x) \in F'$ such that
$\Psi(x)(v) = f_v(x)$ for all $v \in S$, and $\|\Psi(x)\| \leq \|\Phi\|$.
The map $x \mapsto \Psi(x)$ is measurable since the measurable structure
on $F'$ is generated by the linear functionals given by $v \in S$. So
$\Psi \in L^\infty(X, F')$, and regarding the latter as embedded in
$L^\infty(X, F)'$ we have $\Psi(v\cdot 1_X) = \Phi(v\cdot 1_X)$ for all
$v \in S$. But the elements $v\cdot 1_X$ generate
$L^\infty(X, F)$ as a Banach $L^\infty(X)$-module, so this implies that
$\Phi = \Psi$. We conclude that $L^\infty(X, F)' = L^\infty(X, F')$.\hfill\hal
\bigskip

\noindent {\bf Theorem 47.} {\it $\X(F) \cong L^\infty(F, F)$.}
\medskip

Fix $v \in F_0$, and for $t \in \R$ let $\alpha_{tv}: L^\infty(F) \to
L^\infty(F)$ be translation by $tv$. Then $(\alpha_{tv})$ is a strongly
continuous one-parameter group of automorphisms of $L^\infty(F)$, so as in
the proof of Theorem 38 its generator $\delta_v$ is a metric derivation when
restricted to $\Li(F)$.
\medskip

For any $f \in \Li(F)$, we have $\|f - \alpha_{tv}(f)\| \leq L(f)\|tv\|$,
so $\|\delta_v\| \leq \|v\|$. Conversely, find $\phi \in F'$ such that
$\|\phi\| = 1$ and $\phi(v) = \|v\|$, and let $f = (\phi \wedge C)
\vee(-C)$ for some $C > 0$. Then $f \in \Li(F)$ and $L(f) = 1$, but
$(\delta_v f)(w) = \|v\|$ if $-C < \phi(w) < C$. Taking $C \to \infty$, 
we conclude that $|\delta_v| = \|v\|$ almost everywhere.
\medskip

Thus if $A_1, \ldots, A_n \subset F$ are disjoint and
$v_1, \ldots, v_n \in F_0$, then $\sum \chi_{A_i}\delta_{v_i}$ is a metric
derivation and
$$|\sum \chi_{A_i}\delta_{v_i}| = \|v_i\|$$
on $A_i$. Taking completions, this shows that $L^\infty(F, F)$ naturally
isometrically embeds in $\X(F)$.
\medskip

For the converse, we apply Theorem 26 with $E = L^\infty(F, F)$ regarded
as a subset of $\X(F)$. Reflexivity follows from the lemma. To verify the
second hypothesis of Theorem 26, by Theorem 16
it suffices to consider functions in $\Li(F)$ of the form
$f = f_0(\Phi_1, \ldots, \Phi_n)$ for $f_0$ a bounded $C^\infty$ function
on $\R^n$ and $\Phi_1, \ldots, \Phi_n \in F'$. Projecting $F$ onto its
quotient by the intersection of the kernels of $\Phi_1, \ldots, \Phi_n$ and
applying the reasoning in Theorem 44 shows that for such $f$ we have
$|d_E f| = |df|$ as needed.\hfill\hal
\bigskip

\noindent {\bf H. Weiner space}
\bigskip

Consider the Wiener space of continuous functions
$f: [0,\infty) \to \R$ such that $f(0) = 0$. According to ([8], Chapter 3)
there is a measurable isomorphism between this space, equipped with
Weiner measure, and the space $\R^\N =$ the product of countably many
copies of $\R$, giving each factor normalized Gaussian measure. The structure
of the Ornstein-Uhlenbeck operator and the Gross-Sobolev derivative are
more transparent in the $\R^\N$ picture, so we work there. In fact, the
techniques we have introduced in preceding sections suffice to compute
$\X(\R^\N)$ with little further effort.
\medskip

The metric on $\R^\N$ is defined by $\rho(a, b) = (\sum |a_n - b_n|^2)^{1/2}$,
where $a = (a_n)$ and $b = (b_n)$. This metric has infinite distances,
and any ball of finite radius has measure zero, just as in the example
mentioned following Theorem 31. Let $l^2(L^\infty(\R^\N))$ denote the set of
weakly 2-summable sequences $(f_n) \subset L^\infty(\R^\N)$.
\bigskip

\noindent {\bf Theorem 48.} {\it $\X(\R^\N) \cong l^2(L^\infty(\R^\N))$.
This is naturally identified with the space of bounded measurable sections
of the Hilbert bundle $\R^\N\times l^2(\N)$ over $\R^\N$.}
\medskip

\noindent {\it Proof.} For any 2-summable sequence
$(f_n) \in l^2(L^\infty(\R^\N))$ the series
$\sum_0^\infty f_n\partial/\partial x_n$ defines a metric derivation. This is
shown by an argument similar to the one used in the proof of Theorem 44.
In the language of [49], $l^2(L^\infty(\R^\N))$ is a self-dual Hilbert
module over $L^\infty(\R^\N)$. It is therefore reflexive, and
the second hypothesis of Theorem 26 can be verified by the technique of
projection onto finitely many factors also used in the proof of Theorem 44.
\medskip

The realization of $l^2(L^\infty(\R^\N))$ as sections of $\R^\N \times
l^2(\N)$ is given by identifying the 2-summable sequence $(f_n)
\in l^2(L^\infty(\R^\N))$ with the section $x \mapsto
(f_n(x)) \in l^2(\N)$.\hfill\hal
\bigskip

Comparison with [8] and [66] shows that our exterior derivative on
$\R^\N$ is precisely the Gross-Sobolev derivative.
\bigskip
\bigskip

\noindent {\underbar{\bf 6. Dirichlet spaces.}}
\bigskip

\noindent {\bf A. Intrinsic metrics}
\bigskip

In this section we relate our construction to two themes in the study
of Dirichlet spaces. The first is the intrinsic metric associated to
any Dirichlet form. This has been used in several places ([5], [6],
[7], [41], [60-63]), and its geometric aspect was specifically considered in
[60]. The second theme is the existence of a first-order differential
calculus associated to certain Dirichlet spaces; this was hinted at
in [42] and [7] and thoroughly treated in [54] and [55].
\medskip

For basic material on Dirichlet forms we refer the reader to the classic
texts [26] and [59] and the more recent books [8] and [44].
\medskip

The intrinsic metric is most elegantly treated in the setting described in
the next definition. We will connect the structure described here with
traditional Dirichlet forms in Theorem 57.
\bigskip

\noindent {\bf Definition 49.} Let $X$ be a $\sigma$-finite measure space.
An {\it $L^\infty$-diffusion form} is a map
$\Gamma: \D\times\D \to L^\infty(X)$ where $\D$ is a weak*-dense, unital
subalgebra of $L^\infty(X)$, such that
\medskip

{\narrower{
\noindent (a). $\Gamma$ is bilinear, symmetric, and positive (i.e.\
$\Gamma(f,f) \geq 0$ for all $f \in \D$);
\medskip

\noindent (b). $\Gamma(fg, h) = f\Gamma(g,h) + g\Gamma(f,h)$
for all $f,g,h \in \D$ (that is --- by symmetry --- $\Gamma$ is a derivation
in either variable); and
\medskip

\noindent (c). $\Gamma$ is closed in the sense that if $(f_i)
\subset \D$, $\|\Gamma(f_i, f_i)\|$ is uniformly bounded, and
$f_i \to f$ weak* in $L^\infty(X)$, then $f \in \D$ and
$\Gamma(f_i, g) \to \Gamma(f, g)$ weak* in $L^\infty(X)$
for all $g \in \D$.

}}\hfill\hal
\bigskip

For example, if $M$ is a Riemannian manifold then $\D = \Li(M)$ and
$\Gamma(f,g) = \nabla f\cdot \nabla g$ describes an $L^\infty$-diffusion form.
(By Theorems 23 and 36 the standard exterior derivative $d$ is a metric
derivation; hence, identifying the tangent and cotangent bundles, so
is the gradient $\nabla$. This implies the desired closure property.)
In Theorem 56 we generalize this example to include all $L^\infty$-diffusion
forms.
\bigskip

\noindent {\bf Lemma 50.} {\it $|\Gamma(f, g)|^2 \leq \Gamma(f,f)
\Gamma(g,g)$ almost everywhere, for all $f, g \in \D$.}
\medskip

\noindent {\it Proof.} Suppose the inequality fails. Then without loss
of generality there exists
$\epsilon > 0$, a positive measure set $A \subset X$, and scalars
$a, b, c \in \R$ such that $a \geq (bc)^{1/2} + \epsilon$ and
$$\Gamma(f,g) \geq a,\qquad b \geq \Gamma(f,f),\qquad  c \geq \Gamma(g,g)$$
on $A$. Choose $h \in L^1(A)$ with
$h \geq 0$ and $\int h = 1$. Then
integrating $\Gamma$ against $h$ gives rise to a positive semidefinite
bilinear form $\la\cdot, \cdot \ra_h$ on $\D$, and we have
$$\la f, f\ra_h^{1/2} \la g, g\ra_h^{1/2} \leq (bc)^{1/2} \leq a - \epsilon
\leq \la f,g\ra_h - \epsilon,$$
contradicting the Cauchy-Schwartz inequality for $\la\cdot,\cdot\ra_h$.
Thus the desired inequality must hold almost everywhere on $X$.\hfill\hal
\bigskip

\noindent {\bf Theorem 51.} {\it Let $\Gamma$ be an $L^\infty$-diffusion
form. Then there is a measurable metric $\rho$ on $X$ such that
$M = (X, \rho)$ satisfies
$\D = \Li(M)$ and $\|\Gamma(f, f)\|^{1/2} = L(f)$ for all
$f \in \D$. $M$ is differentiable in the sense of Definition 30.}
\medskip

\noindent {\it Proof.} For each $g \in \D$ with $\|\Gamma(g, g)\| \leq 1$
let $X_g$ be a copy of $X$; then let $Y = \bigcup X_g$ be their disjoint
union. $L^\infty(Y)$ is an abelian W*-module over $L^\infty(X)$ via the
diagonal embedding of $L^\infty(X)$ into $L^\infty(Y)$.
Define $\delta: \D \to L^\infty(Y)$ by $\delta(f) = \Gamma(f,g)$ on
$X_g$. Then $\delta$ is a W*-derivation and so Theorem 19 implies the
existence of a measurable metric $\rho$ on $X$ such that $M = (X,\rho)$
satisfies $\D = \Li(M)$ and $L(f) = \|\delta f\|$ for all $f \in \D$.
\medskip

By the lemma we have $|\delta f| \leq \Gamma(f,f)^{1/2}$ almost everywhere
on $X$. Conversely, taking $g = f/\|\Gamma(f,f)\|^{1/2}$ we have
$$\|\delta f\| \geq \|\Gamma(f,f)\|/\|\Gamma(f,f)\|^{1/2} =
\|\Gamma(f,f)\|^{1/2}.$$
So $L(f) = \|\delta f\| = \|\Gamma(f,f)\|^{1/2}$.\hfill\hal
\bigskip

Concretely, the measurable metric identified in Theorem 51 is given by
$$\rho(A, B) = \sup\{\rho_f(A, B): f \in \D, \|\Gamma(f,f)\| \leq 1\}.$$
Thus, it is essentially the intrinsic metric mentioned earlier, except that
the latter is a pointwise metric and cannot be defined without regularity
assumptions sufficient to make every $f \in \D$ well-defined
at each point of $X$. This can always be assured by altering the underlying
space via Theorem 20.
\bigskip

\noindent {\bf B. Tangent and cotangent bundles}
\bigskip

We now present Sauvageot's construction of an exterior derivative in the
setting of $L^\infty$-diffusion forms and compare it to our exterior
derivative. The approach of [54] has been
altered slightly here to more clearly display its connection with K\"ahler
differentiation; see e.g.\ ([31], $\S$ 20) or ([33], $\S$ II.8).
\medskip

\noindent {\bf Definition 52.} Let $\Gamma$ be an $L^\infty$-diffusion
form on $L^\infty(X)$. Let $E_0$ be the algebraic tensor product
$E_0 = \D \otimes \D$, equipped with the $L^\infty(X)$-valued inner product
$$\la f_1 \otimes g_1, f_2\otimes g_2\ra_{L^\infty} = f_1f_2\Gamma(g_1, g_2)$$
and regarded as a bimodule over $\D$ with left and right actions given by
$f\cdot (g\otimes h)\cdot k = fg \otimes hk$.
\medskip

Let $I = \{\phi \in E_0: \la \phi, \phi\ra_{L^\infty} = 0\}$ and define $E_1$
to be the sub-bimodule of $E_0/I$ generated by the elements of the form
$f\otimes 1 - 1\otimes f$ for $f \in \D$. A short calculation using the
derivation identity (Definition 49 (b)) shows that $E_1$ is a monomodule,
i.e.\ $f\phi = \phi f$ for all $f \in \D$ and $\phi \in E_1$.
\medskip

Finally, let $E$ be the set of bounded $\D$-module homomorphisms from
$E_1$ to $L^\infty(X)$. Note that $E_1$ naturally embeds in $E$ by identifying
$\phi \in E_1$ with the homomorphism $\psi \mapsto \la \psi,\phi\ra_{L^\infty}$.
Furthermore, $E$ is an $L^\infty(X)$-module with action
$f\cdot \Phi(\phi) = f\Phi(\phi)$ ($f \in L^\infty(X)$, $\Phi \in
E$, $\phi \in E_1$) and is a self-dual Hilbert module by
([49], Theorem 4.2).\hfill\hal
\bigskip

\noindent {\bf Definition 53.} Retain the notation used
in Definition 52. Define $\delta: \D \to E$ by
$$\delta(f) = 1 \otimes f - f \otimes 1.$$
It is easy to check that $\delta$ is a derivation.\hfill\hal
\bigskip

\noindent {\bf Proposition 54.} {\it Let $\Gamma$ be an $L^\infty$-diffusion
form and let $M = (X,\rho)$, $E$, and $\delta$ be the associated structures
identified in Theorem 51 and Definitions 52 and 53. Then $\Gamma(f,g) =
\la \delta f, \delta g\ra_{L^\infty}$ for all $f, g \in \D$, $\delta$ is a metric
derivation, and $\|\delta f\| = L(f)$ for all $f \in \D = \Li(M)$. There
is a nonexpansive
module homomorphism $T: E \to \X(M)$ such that $\delta = T^*\circ d$.}
\medskip

\noindent {\it Proof.} Since $\Gamma(1, f) = 0$ for any $f \in \D$ by the
derivation property, we have
$$\la \delta f, \delta g\ra_{L^\infty} = \la 1\otimes f, 1\otimes g\ra_{L^\infty} =
\Gamma(f,g)$$
for all $f, g \in \D$. It follows from Theorem 51 that $\|\delta f\| =
\|\Gamma(f,f)\|^{1/2} = L(f)$ for all $f \in \D$. To see that $\delta$
is a metric derivation, suppose $f_i \to f$ boundedly weak* in $\Li(M)$
and let $g, h \in \Li(M)$. Then
$$\la \delta f_i, g\otimes h\ra_{L^\infty} =
g\Gamma(f_i, h) \to g\Gamma(f, h) = \la \delta f, g\otimes h\ra_{L^\infty}$$
weak* in $L^\infty(X)$, and this is sufficient to show that $\delta f_i
\to \delta f$
weak* in $E$ ([49], Remark 3.9). So $\delta$ is a metric derivation.
\medskip

$T: E \to \X(M)$ is defined by $(T\phi)(f) = \la \delta f, \phi\ra_{L^\infty}$ for
$\phi \in E$ and $f \in \Li(M)$. For any $f \in \Li(M)$ we then have
$$(T^*\circ d)(f)(\phi) = df(T\phi) = (T\phi)(f) = \la \delta f, \phi \ra_{L^\infty},$$
so that $(T^*\circ d)(f) = \delta f$.\hfill\hal
\bigskip

In most natural examples the map $T$ in Proposition 54 is an isometric
isomorphism, so that $\X(M) \cong \Omega(M)$ is the module of bounded
measurable sections of a bundle of Hilbert spaces (see Corollary 24 and
the comment following Theorem 10). However, this is not always the case.
We now give two examples to illustrate what other behavior can occur.
(Besides these examples, it is also instructive to consider abstract
Weiner spaces, as discussed e.g.\ in [8]. There one has a Hilbert space
$H$, a Banach space $F$, and maps $H \to F$ and $F' \to H$ which
correspond to the maps $T$ and $T^*$ of Proposition 54 acting on the
tangent and cotangent spaces at a single point.)
\medskip

Let $M_1 = [0,1]^2$ be the unit square with Euclidean metric and
Lebesgue measure. Let
$S = M_1 \cap \Q^2$ be the set of points with rational coordinates.
For each $s = (x,y) \in S^2$ let $O_s$ be an open subset of $M_1$ which
contains the line segment joining $x$ and $y$; since $S^2$ is countable we
can arrange that the total measure of $A = \bigcup_{s \in S^2} O_s$ is
less than one (indeed, arbitrarily close to zero). Let $h \in L^\infty(M_1)$
satisfy $0 \leq h \leq 1$ and $h|_A = 1$, and define
$$\Gamma_1(f,g) = h\nabla f\cdot \nabla g$$
for $f, g \in \Li(M_1)$. Let $E$ and $\delta$ be the associated Hilbert
module and metric derivation identified in Definitions 52 and 53, and let
$d: \Li(M_1) \to \Omega(M_1)$ be the exterior derivative.
\medskip

Let $M_2 = [0,1]^2$ be the unit square with $l^1$ metric, i.e.\
$$\rho(x,y) = |x_1 - y_1| + |x_2 - y_2|$$
for $x,y \in M_2$. Let $B_h = (\R\times \Q) \cap M_2$, $B_v =
(\Q \times \R)\cap M_2$, and $B = B_h \cup B_v$. Then $B$
is a countable union of horizontal and vertical line segments, so we can
find a probability measure $\mu$ on $M_2$ which is supported on $B$ and
whose restriction to any line segment is
a positive multiple of Lebesgue measure. For $f, g \in \Li(M_2)$ define
$$\Gamma_2(f,g) = \cases{
\big({{\partial f}\over{\partial x_1}}\big)
\big({{\partial g}\over{\partial x_1}}\big)&on $B_h$\cr
\big({{\partial f}\over{\partial x_2}}\big)
\big({{\partial g}\over{\partial x_2}}\big)&on $B_v$\cr
0&elsewhere\cr}.$$
(Note that $B_h \cap B_v$ has measure zero, so the definition of
$\Gamma_2$ is consistent.)
\bigskip

\noindent {\bf Theorem 55.} {\it $\Gamma_1$ and $\Gamma_2$ are
$L^\infty$-diffusion forms and $M_1$ and $M_2$ are the associated metric
spaces. For $f \in L^\infty(M_1)$ and $g \in \Li(M_1)$ we have
$|f dg| = f|\nabla g|$ regardless of $h$ but
$|f\delta g| = f h^{1/2}|\nabla g|$.}
\medskip

\noindent {\it Proof.} In both cases the metric of Theorem 51 agrees with the
given metric on $[0,1]^2 \cap \Q^2$, hence $L(f) = \|\Gamma(f,f)\|^{1/2}$
for all $f \in \D$. The rest of the proof is routine.\hfill\hal
\bigskip

Thus, taking $f \in L^\infty(M_1)$ supported on $\{x: h(x) \leq 1/2\}$, say,
and $g \in \Li(M_1)$ with $|\nabla g| = 1$ almost everywhere, we see
that for $M_1$ the map $T$ of Proposition 54 is not isometric; and if
$h$ is zero on a positive measure set $T^*$ will have a nonzero kernel.
Whereas $\X(M_2)$ and $\Omega(M_2)$ are not even Hilbert modules, so $T$
certainly cannot be isometric in this case.
\medskip

$M_2$ can be viewed as a (non-smooth) sub-Riemannian metric, as can $M_1$
in the case $h = \chi_A$.
\medskip

The next result, an easy corollary of Proposition 54, characterizes
$L^\infty$-diffusions in metric terms.
\bigskip

\noindent {\bf Theorem 56.} {\it Let $M$ be a measurable metric space,
let $E$ be a self-dual Hilbert module over $L^\infty(M)$, and
let $\delta: M \to E$ be a metric derivation which satisfies $\|\delta f\|
= L(f)$ for all $f \in \Li(M)$. Then $\D = \Li(M)$ and
$\Gamma(f,g) = \la \delta f, \delta g\ra_{L^\infty}$
define an $L^\infty$-diffusion form on $L^\infty(M)$. Conversely, every
$L^\infty$-diffusion form arises in this manner.}
\medskip

\noindent {\it Proof.} It is easy to check that $\Gamma$ is an
$L^\infty$-diffusion form. The closure property of Definition 49 (c)
follows from the hypothesis $\|\delta f\| = L(f)$ (so that bounded
$L^\infty$-weak* convergence in $\D$ is weak*-convergence in $\Li(M)$)
plus the fact that $\delta$ is weak*-continuous.
\medskip

The converse statement is Proposition 54.\hfill\hal
\bigskip

\noindent {\bf C. Dirichlet forms}
\bigskip

We adopt the conventions of [8]. Thus, if $X$ is a $\sigma$-finite measure
space then a {\it Dirichlet form} on $L^2(X)$ is a positive, symmetric,
closed bilinear map $\E: D(\E)\times D(\E) \to \R$, where $D(\E)$ is a dense
subspace of $L^2(X)$, such that $f\in D(\E)$ implies $f \wedge 1 \in D(\E)$
and
$$\E(f\wedge 1, f\wedge 1) \leq \E(f,f).$$
Associated to any Dirichlet form $\E$ there is a {\it sub-Markovian
symmetric semi-group} $(P^{(2)}_t)$ ($t\geq 0$), which is a norm continuous
semigroup of self-adjoint contractions on $L^2(X)$ that satisfy
$0 \leq P^{(2)}_tf \leq 1$
whenever $0 \leq f \leq 1$. The infinitesimal generator
$A^{(2)}$ of the semigroup corresponding to the Dirichlet form $\E$ satisfies
$D((-A^{(2)})^{1/2}) = D(\E)$ and $\E(f, g) = -\la A^{(2)}f, g\ra$ for all
$f\in D(A^{(2)})$ and $g \in D(\E)$.
\medskip

The restrictions of the operators $P^{(2)}_t$ to $L^1(X)\cap L^2(X)$ extend to
a norm continuous contraction semigroup $P^{(1)}_t$ on $L^1(X)$ whose
infinitesimal generator is denoted $A^{(1)}$. The
intersection of $D(A^{(2)})$ with $D(A^{(1)})$ is a core for both $A^{(1)}$
and $A^{(2)}$, and on this intersection the two operators agree.
\medskip

$\E$ is a {\it diffusion form} if its jump and killing parts are zero
([59], [26]). If the form is regular, this is equivalent to strong locality
(follows from Theorem 11.10 of [59]; see also section 4.1.i of [61]).
\medskip

The Dirichlet form $\E$ is said to {\it admit a carr\'e du champ} if
$D(A^{(1)}) \cap L^\infty(X)$ is an algebra. (See [8] for several
equivalent conditions.) In this case there is a positive, symmetric,
bilinear map $\Gamma^{(2)}: D(\E)\times D(\E) \to L^1(X)$ which satisfies
$$\E(fh,g) + \E(gh, f) - \E(h, fg) = \int h\Gamma^{(2)}(f, g)$$
for all $f, g, h \in D(\E) \cap L^\infty(X)$. If $\E$ is a diffusion form,
this can be simplified to the equality
$$2\E(f,g) = \int \Gamma^{(2)}(f,g)$$
for all $f,g \in D(\E)$.
\medskip

We now show how to derive $L^\infty$-diffusion forms from
diffusion forms which admit a {\it carr\'e du champ}. Morally, one should
be able to go in the reverse direction as well by defining
$\E(f,g) = \int \Gamma(f,g)/2$ when an $L^\infty$-diffusion form $\Gamma$
is given. Two difficulties arise, however. First, there is the question of
whether $\Gamma(f,g)$ is integrable for sufficiently many $f$ and $g$ (a
problem that can presumably be addressed by taking some care in choosing
$\mu$ from its measure class); more seriously, it seems that the closure
condition on $\Gamma$ is weaker than the corresponding closure condition
on $\E$. So there may be cases where one cannot pass in this direction,
though I do not know of any examples.
\bigskip

\noindent {\bf Theorem 57.} {\it Let $X$ be a $\sigma$-finite measure space
and let $\E$ be a diffusion form on $L^2(X)$ which admits a
{\it carr\'e du champ}. Then there is a unique $L^\infty$-diffusion
form $\Gamma$ which agrees with $\Gamma^{(2)}$ on a common core.}
\medskip

\noindent {\it Proof.}
Define a weak*-continuous contraction semigroup $P_t^{(\infty)}$
on $L^\infty(X)$ by letting $P_t^{(\infty)}$ be the adjoint of $P_t^{(1)}$.
For any $f, g \in L^1(X) \cap L^\infty(X)$ we have $\int P_t^{(2)}(f)g =
\int fP_t^{(2)}(g)$ since $P_t^{(2)}$ is self-adjoint; and since $P_t^{(1)}$
agrees with $P_t^{(2)}$ on $L^1(X) \cap L^\infty(X)$
this implies that $P_t^{(\infty)} =
P_t^{(1)}$ on $L^1(X) \cap L^\infty(X)$. Thus, the infinitesimal generator
$A^{(\infty)}$ of $(P_t^{(\infty)})$ agrees with $A^{(1)}$ on
$S = D(A^{(1)}) \cap L^\infty(X)$.
\medskip

For $f, g \in S$ define
$$\Gamma(f,g) = -(A^{(\infty)}(fg) - fA^{(\infty)}(g) - gA^{(\infty)}(f)).$$
This makes sense because $S$ is an algebra by ([8], Theorem 4.2.1), and
$\Gamma = \Gamma^{(2)}$ on $S$ by ([8], Theorem 4.2.2).
The derivation property follows from the fact that $\E$ is a
diffusion form by ([42], $\S$ 1.5); the remaining issue is closability
of $\Gamma$.
\medskip

For $A, B \subset X$ define
$$\rho(A, B) = \sup\{\rho_f(A, B): f \in S\hbox{ and }\Gamma(f,f) \leq 1\}.$$
By Theorem 16, $S$ is weak*-dense in $\Li(M)$ where $M = (X,\rho)$. We must
define $\Gamma(f,g)$ for all $f,g \in \Li(M)$.
\medskip

First, for $f \in \Li(M)$ and $g \in S$, set
$$\Gamma(f,g) = \lim \Gamma(f_i, g)$$
where $(f_i) \subset S$ and $f_i \to f$ boundedly weak* in $\Li(M)$. To
see that this definition makes sense, we must show that $f_i \to 0$
implies $\Gamma(f_i, g) \to 0$. Suppose $f_i \to 0$ weak* and
let $h \in S$; then
$$\eqalign{\int \Gamma(f_i, g)h
&= \int -A^{(\infty)}(f_ig)h + f_iA^{(\infty)}(g)h + gA^{(\infty)}(f_i)h\cr
&= \int -f_igA^{(1)}(h) + f_iA^{(\infty)}(g)h + f_iA^{(1)}(gh)\cr
&\to 0,\cr}$$
as desired.
Thus $\Gamma(f,g)$ is well-defined when $f \in \Li(M)$ and $g \in S$.
\medskip

To define $\Gamma(f,g)$ for any $f,g \in \Li(M)$ we use the key trick in
the proof of ([49], Theorem 3.2). Fix $f, g \in \Li(M)$. For any
$h \in L^1(X)$, $h \geq 0$,
let $H_h$ be the pre-Hilbert space consisting of the set $S$ with the
pseudonorm $\|k\|^2 = \int \Gamma(k,k)h$.
Then the maps $\phi_f: k \mapsto \int \Gamma(f,k)h$ and $\phi_g:
k \mapsto \int \Gamma(g,k)h$ are bounded linear functionals on $H_h$,
hence both are represented by elements of the Hilbert space completion
of $H_h$. Define $\Gamma(f,g)_h$ to be the inner product of $\phi_f$
and $\phi_g$, and observe that $|\Gamma(f,g)_h| \leq L(f)L(g)\|h\|_1$.
We can then define $\Gamma(f,g)_h$ for all $h \in L^1(X)$ by linearity
and let $\Gamma(f,g) \in L^\infty(X)$ be defined by
$$\int \Gamma(f,g)h = \Gamma(f,g)_h.$$
\medskip

$\Gamma$ is now closed because $f_i \to f$ boundedly weak* in $\Li(M)$
implies $\phi_{f_i} \to \phi_f$ for all $h \in S$, hence
$\Gamma(f_i, g)_h \to \Gamma(f,g)_h$. The other desired properties of
$\Gamma$ hold by continuity.\hfill\hal
\bigskip
\bigskip

[1] M.\ T.\ Barlow and R.\ F.\ Bass, The construction of Brownian motion
on the Sierpinski carpet, {\it Ann.\ Inst.\ Henri Poincar\'e \bf 25}
(1989), 225-257.
\medskip

[2] M.\ T.\ Barlow and E.\ A.\ Perkins, Brownian motion on the Sierpinski
gasket, {\it Probab.\ Th.\ Rel.\ Fields \bf 79} (1988), 543-623.
\medskip

[3] A.\ Bella\"iche, The tangent space in sub-Riemannian geometry, in [4],
1-78.
\medskip

[4] A.\ Bella\"iche and J.-J.\ Risler, eds., {\it Sub-Riemannian Geometry},
Birkh\"auser (1996).
\medskip

[5] A.\ D.\ Bendikov, Remarks concerning the analysis on local Dirichlet
spaces, in {\it Dirichlet Forms and Stochastic Processes (Beijing, 1993)},
de Gruyter (1995), 55-64.
\medskip

[6] M.\ Biroli and U.\ Mosco, Formes de Dirichlet et estimations
structurelles dans les milieux discontinus, {\it C.\ R.\ Acad.\ Sci.\
Paris \bf 313} (1991), 593-598.
\medskip

[7] ---------, A Saint-Venant type principle for Dirichlet forms on
discontinuous media, {\it Ann.\ Mat.\ Pura Appl.\ \bf 169} (1995), 125-181.
\medskip

[8] N.\ Bouleau and F.\ Hirsch, {\it Dirichlet Forms and Analysis on
Wiener Space}, de Gruyter (1991).
\medskip

[9] B.\ H.\ Bowdich, Notes on locally CAT(1) spaces, in {\it Geometric
Group Theory} (ed.\ R.\ Charney, M.\ Davis, and M.\ Shapiro), de Gruyter
(1995), 1-48.
\medskip

[10] O.\ Bratteli and D.\ W.\ Robinson, {\it Operator Algebras and
Quantum Statistical Mechanics 1} (second edition), Springer-Verlag (1987).
\medskip

[11] H.\ Cartan and S.\ Eilenberg, {\it Homological Algebra}, Princeton (1956).
\medskip

[12] A. Connes, C* alg\'ebres et g\'eom\'etrie diff\'erentielle,
{\it C.\ R.\ Acad.\ Sc.\ Paris Ser.\ A \bf 290} (1980), 599-604.
\medskip

[13] ---------, {\it Noncommutative Geometry}, Academic Press (1994).
\medskip

[14] J.\ B.\ Conway, {\it A Course in Functional Analysis}, Springer GTM
{\bf 96} (1985).
\medskip

[15] E.\ B.\ Davies, Large deviations for heat kernels on graphs,
{\it J.\ London Math.\ Soc.\ \bf 47} (1993), 65-72.
\medskip

[16] ---------, Analysis on graphs and noncommutative geometry,
{\it J.\ Funct.\ Anal.\ \bf 111} (1993), 398-430.
\medskip

[17] G.\ De Cecco and G.\ Palmieri, LIP manifolds: from metric to Finslerian
structure, {\it Math.\ Zeit.\ \bf 218} (1995), 223-237.
\medskip

[18] M.\ Dubois-Violette, D\'erivations et calcul diff\'erentiel non
commutatif, {\it C.\ R.\ Acad.\ Sci.\ Paris \bf 307}, (1988), 403-408.
\medskip

[19] M.\ J.\ Dupr\'e and R.\ M.\ Gillette, {\it Banach bundles, Banach
modules and automorphisms of C*-algebras}, Pitman Research Notes in
Mathematics {\bf 92} (1983).
\medskip

[20] J.\ Eells, Jr., On submanifolds of certain function spaces,
{\it Proc.\ Nat.\ Acad.\ Sci.\ \bf 45} (1959), 1520-1522.
\medskip

[21] E.\ G.\ Effros and Z.-J.\ Ruan, Representations of operator bimodules
and their applications, {\it J.\ Op.\ Thy.\ \bf 19} (1988), 137-157.
\medskip

[22] H.\ Federer, {\it Geometric Measure Theory}, Springer-Verlag (1969).
\medskip

[23] J.\ M.\ G.\ Fell, {\it Induced Representations and Banach
*-Algebraic Bundles}, Springer-Verlag LNM {\bf 582} (1977).
\medskip

[24] R.\ Feres, personal communication.
\medskip

[25] G.\ B.\ Folland, {\it Real Analysis}, Wiley-Interscience (1984).
\medskip

[26] M.\ Fukushima, {\it Dirichlet Forms and Markov Processes},
North-Holland (1980).
\medskip

[27] V.\ M.\ Gol'dshtein, V.\ I.\ Kuz'minov, and I.\ A.\ Shvedov,
Differential forms on Lipschitz manifolds, {\it Siberian Math.\ J.\ \bf 23}
(1982), 151-161.
\medskip

[28] K.\ Goodearl, personal communication.
\medskip

[29] M.\ Gromov, {\it Structures M\'etriques Pour les Vari\'et\'es
riemanniennes}, Cedic/
Fernand Nathan (1981).
\medskip

[30] ---------, Hyperbolic groups, in {\it Essays in Group Theory} (ed.\
S.\ M.\ Gersten), Springer-Verlag (1988), 75-263.
\medskip

[31] A.\ Grothendieck, {\it Fondements de la G\'eom\'etric Alg\'ebrique}
volume IV, {\it Publ.\ Math.\ IHES \bf 20} (1962).
\medskip

[32] J.\ Harrison, Stokes' theorem for nonsmooth chains, {\it Bull.\ Amer.\
Math.\ Soc.\ \bf 29} (1993), 235-242.
\medskip

[33] R.\ Hartshorne, {\it Algebraic Geometry}, Springer-Verlag GTM {\bf 52}
(1977).
\medskip

[34] P.\ E.\ T.\ Jorgensen and S.\ Pedersen, Harmonic analysis of fractal
measures induced by representations of a certain C*-algebra,
{\it Bull.\ Amer.\ Math.\ Soc.\ \bf 29} (1993), 228-234.
\medskip

[35] R.\ V.\ Kadison and J.\ R.\ Ringrose, {\it Fundamentals of the Theory
of Operator Algebras} (volume II), Academic Press (1986).
\medskip

[36] J.\ Kigami, A harmonic calculus on the Sierpinski spaces,
{\it Japan J.\ Appl.\ Math.\ \bf 6} (1989), 259-290.
\medskip

[37] ---------, Harmonic calculus on P.C.F.\ self-similar sets,
{\it Trans.\ Amer.\ Math.\ Soc.\ \bf 335} (1993), 721-755.
\medskip

[38] A.\ Kor\`anyi, M.\ A.\ Picardello, and M.\ H.\ Taibleson,
Hardy spaces on non-homogeneous trees, {\it Sympos.\ Math.\ \bf 29}
(1987), 205-265.
\medskip

[39] S.\ Kusuoka, Dirichlet forms on fractals and products of random matrices,
{\it Publ.\ RIMS, Kyoto Univ.\ \bf 25} (1989), 659-680.
\medskip

[40] G.\ Landi, Nguyen Ai V., K.\ C.\ Wali, Gravity and electromagnetism
in noncommutative geometry, {\it Phys. Lett. B \bf 326} (1994), 45-50.
\medskip

[41] M.\ Lapidus, Analysis on fractals, Laplacians on self-similar sets,
noncommutative geometry and spectral dimensions, {\it Topol.\ Methods
Nonlinear Anal.\ \bf 4} (1994), 137-195.
\medskip

[42] Y.\ Lejan, Mesures associ\'ees a une forme de Dirichlet. Applications,
{\it Bull.\ Soc.\ Math.\ France \bf 106} (1978), 61-112.
\medskip

[43] J.\ Luukkainen and J.\ V\"ais\"al\"a, Elements of Lipschitz topology,
{\it Ann.\ Acad.\ Sci.\ Fenn.\ Ser.\ A  I \bf 3} (1977), 85-122.
\medskip

[44] Z.-M.\ Ma and M.\ R\"ockner, {\it Introduction to the Theory of
(Non-Symmetric) Dirichlet Forms}, Springer-Verlag (1991).
\medskip

[45] J.\ Madore, {\it An Introduction to Noncommutative Differential
Geometry and its Physical Applications}, LMS Lecture Note Series {\bf 206}
(1995).
\medskip

[46] P.\ Mattila, {\it Geometry of Sets and Measures in Euclidean Spaces},
Cambridge (1995).
\medskip

[47] J.\ Milson and B.\ Zombro, A K\"ahler structure on the moduli space of
isometric maps of a circle into Euclidean space, {\it Invent.\ Math.\
\bf 123} (1996), 35-59.
\medskip

[48] P.\ Pansu, M\'etriques de Carnot-Carath\'eodory et quasiisom\'etries
des espaces

\noindent sym\'etriques de rang un, {\it Ann.\ Math.\ \bf 129} (1989), 1-60.
\medskip

[49] W.\ L.\ Paschke, Inner product modules over B*-algebras, {\it Trans.\
Amer.\ Math.\ Soc.\ \bf 182} (1973), 443-468.
\medskip

[50] P.\ Petersen, {\it Riemannian Geometry}, Springer-Verlag (1997).
\medskip

[51] N.\ C.\ Phillips and N.\ Weaver, Modules with norms which take values
in a C*-algebra, {\it Pac.\ J.\ Math.}, to appear.
\medskip

[52] J.\ Rosenberg, Applications of analysis on Lipschitz manifolds,
in {\it Miniconferences on Harmonic Analysis and Operator Algebras
(Canberra, 1987)}, {\it Proc.\ Centre Math.\ Anal.\ Austral.\ Nat.\ Univ.
\bf 16}, Austral.\ Nat.\ Univ.\ (1988).
\medskip

[53] C.\ Sabot, Existence et unicit\'e de la diffusion sur un ensemble fractal,
{\it C.\ R.\ Acad.\ Sci.\ Paris Ser.\ I \bf 321} (1995), 1053-1059.
\medskip

[54] J.-L.\ Sauvageot, Tangent bimodule and locality for dissipative
operators on C*-algebras, in {\it Quantum Probability and Applications IV},
Springer LNM {\bf 1396} (1989), 322-338.
\medskip

[55] ---------, Quantum Dirichlet forms, differential calculus and
semigroups, in {\it Quantum Probability and Applications V}, Springer LNM
{\bf 1442} (1990), 334-346.
\medskip

[56] S.\ Semmes, Finding curves on general spaces through quantitative
topology, with applications to Sobolev and Poincare inequalities,
{\it Sepecta Math.\ \bf 2} (1996), 155-295.
\medskip

[57] R.\ Sikorski, Abstract covariant derivative, {\it Colloq.\ Math.\
\bf 18} (1967), 251-272.
\medskip

[58] ---------, Differential modules, {\it Colloq.\ Math.\ \bf 24} (1971/72),
45-79.
\medskip

[59] M.\ L.\ Silverstein, {\it Symmetric Markov Processes}, Springer LNM
{\bf 426} (1974).
\medskip

[60] K.-T.\ Sturm, On the geometry defined by Dirichlet forms, in {\it
Seminar on Stochastic Analysis, Random Fields and Applications (Ascona, 1993)},
{\it Progr.\ Probab.\ \bf 36} (1995), 231-242.
\medskip

[61] ---------, Analysis on local Dirichlet spaces I. Recurrence,
conservativeness and $L^p$-Liouville properties, {\it J.\ Reine Agnew.\
Math.\ \bf 456} (1994), 173-196.
\medskip

[62] ---------, Analysis on local Dirichlet spaces II. Upper Gaussian
estimates for the fundamental solutions of parabolic equations,
{\it Osaka J.\ Math.\ \bf 32} (1995), 275-312.
\medskip

[63] ---------, Analysis on local Dirichlet spaces III. The parabolic
Harnack inequality, {\it J.\ Math.\ Pures Appl.\ \bf 75} (1996), 273-297.
\medskip

[64] M.\ H.\ Taibleson, Hardy spaces of harmonic functions on homogeneous
isotropic trees, {\it Math.\ Nachr.\ \bf 133} (1987), 273-288.
\medskip

[65] A.\ Takahashi, {\it Fields of Hilbert Modules}, Dissertation,
Tulane University (1971).
\medskip

[66] A.\ S.\ \"Ust\"unel, {\it An Introduction to Analysis on Wiener space},
Springer-Verlag LNM {\bf 1610} (1995).
\medskip

[67] N.\ Weaver, Order completeness in Lipschitz algebras, {\it J.\ Funct.\
Anal.\ \bf 130} (1995), 118-130.
\medskip

[68] ---------, Nonatomic Lipschitz spaces, {\it Studia Math.\ \bf 115}
(1995), 277-289.
\medskip

[69] ---------, Lipschitz algebras and derivations of von Neumann algebras,
{\it J.\ Funct.\ Anal.\ \bf 139} (1996), 261-300.
\bigskip
\bigskip

\noindent Mathematics Dept.

\noindent Washington University

\noindent St.\ Louis, MO 63130 USA

\noindent nweaver@math.wustl.edu
\end

Let $E$ be the Hilbert module associated to $\Gamma$ in Definition 52, and
define $\la \phi, \psi\ra = \int \la \phi, \psi\ra_{L^\infty} d\mu$
for $\phi, \psi \in E$ such that $\int |\phi|^2d\mu, \int |\psi|^2d\mu
< \infty$. This makes the set of such elements of $E$ into a pre-Hilbert
space $H$, and we have $\E(f, g) = \la df, dg\ra$ for all $f, g \in S$.
Thus our hypothesis says that $(df_n)$ is Cauchy in $H$, and we must
show that $df_n \to 0$. Without loss of generality suppose $\|df_n\|_2
\leq 1$ for all $n$.
\medskip

Let $N \in \N$ and $\epsilon > 0$. By passing to a subsequence, we may
assume that $\|f_n\|_2^2 \leq \epsilon/2^n$. Thus
$$\mu(\{x: |f_n(x)| \geq 1\}) \leq \epsilon/2^n,$$
and so
$$A_1 = \{x: |f_n(x)| \geq 1\hbox{ for some }n\}$$
has measure at most $\epsilon$. Similarly, we may assume that
$$A_2 = \{x: |df_n(x)| \geq N\hbox{ for some }n\}$$
has measure at most $1/N$. Now for any $A \subset X$
the map $\phi \mapsto \chi_{B}\phi$ defines a bounded
operator on the pre-Hilbert space $H$; since $\epsilon$ can be arbitrarily
small and $N$ arbitrarily large, it will
suffice to show that $\chi_B df_n \to 0$ for any $B \subset X$ with
finite measure on which $|f_n| \leq 1$ and $|df_n| \leq N$ for all $n$.
\medskip

Let $g_n = f_n|_B$.